\documentclass[twoside,10pt]{article}
\usepackage{mathrsfs}
\usepackage{amssymb}
\usepackage{amsmath}
\usepackage{amsfonts}
\usepackage{amssymb, amsmath,cite,color}
\usepackage[top=1in,bottom=1in,left=0.80in,right=0.80in]{geometry}

\numberwithin{equation}{section} \numberwithin{theorem}{section}
\numberwithin{definition}{section} \numberwithin{lemma}{section}
\numberwithin{corollary}{section} \numberwithin{remark}{section}
\numberwithin{example}{section} \numberwithin{Claim}{section}

\begin{document}
\baselineskip 20pt
\title{{\bf  Monotone iterative schemes for positive
solutions of a fractional differential system with integral boundary
conditions on an infinite interval
 }}
    \author{Yaohong Li$^{a}$, Wei Cheng$^{b}$,
    \ Jiafa Xu$^b$\footnote{Corresponding author: Jiafa Xu, School of Mathematical Sciences, Chongqing Normal University, Chongqing 401331, People's Republic of China.}
    \footnote{E-mail addresses: liz.zhanghy@163.com(Y.Li), 1375415619@qq.com(W.
Cheng), 20150028@cqnu.edu.cn(J.Xu).}\\[5pt]
 $^a$ \small School of Mathematics and Statistics, Suzhou
University,\\
\small Suzhou 234000, Anhui, P. R. China.\\[5pt]
  $^b$ \small School of Mathematical Sciences, Chongqing Normal University, \\ \small  Chongqing 401331, P. R. China.}
\date{}
 \maketitle
\begin{abstract} In this paper,
using the monotone iterative technique and the Banach contraction
mapping principle, we study a class of fractional differential
system with integral boundary on an infinite interval. Some explicit
monotone iterative schemes for approximating the extreme positive
solutions and the unique positive solution are constructed.

\textbf{Keywords:} monotone iterative technique; iterative schemes;
fractional differential system; infinite interval.

\textbf{MR(2010) Subject Classification:} 34A05; 34B18; 26A33.

\end{abstract}

\section{Introduction}

The purpose of this paper is to study monotone iterative schemes of
positive solutions for the following fractional differential system
with integral boundary conditions
\begin{equation}\label{fds} \left\{
      \begin{array}{ll}
       D^{\alpha_{1}}u(t)+f_{1}(t,u(t),v(t),D^{\alpha_{1}-1}u(t),D^{\alpha_{2}-1}v(t))=0,\
       n_{1}-1<\alpha_{1}\leq n_{1},\\
        D^{\alpha_{2}}v(t)+f_{2}(t,u(t),v(t),D^{\alpha_{1}-1}u(t),D^{\alpha_{2}-1}v(t))=0,\
        n_{2}-1<\alpha_{2}\leq n_{2},\\
       u(0)=u'(0)=\cdots=u^{(n_{1}-2)}(0)=0,D^{\alpha_{1}-1}u(+\infty)=
       \displaystyle\int_{0}^{+\infty}h_{1}(t)u(t)\mbox{d}t,\\
        v(0)=v'(0)=\cdots=v^{(n_{2}-2)}(0)=0,D^{\alpha_{2}-1}v(+\infty)=
       \displaystyle\int_{0}^{+\infty}h_{2}(t)v(t)\mbox{d}t,\\
     \end{array}
   \right.\end{equation}
  where $t\in J=[0,+\infty),f_{i}\in
C(J\times\mathbb{R}\times\mathbb{R}\times\mathbb{R}\times\mathbb{R},J),n_{i}\in
N^{+},h_{i}(t)\in L[0,+\infty), D^{\alpha_{i}}$ are the standard
Riemann-Liouville fractional derivative of order $\alpha_{i},i=1,2$.
Here we emphasize that the nonlinearity terms $f_{i}$ rely on the
lower-order fractional derivative of multiple unknown functions and
the fractional infinite boundary value rely on the infinite integral
of unknown functions.

 In recent decades, there has been a rapid growth in the number of fractional calculus
 from both theoretical and applied perspectives, more detailed description of the subject can be found
 in the books \cite{bk1,bk2,bk3,bk4}.
We note that most of the current results on the
 existence of fractional differential equations are focused on the finite
interval, see
\cite{c5,c6,c7,c8,c9,c10,c11,c12,c13,c14,c15,c16,c17,c18,c19,c20,c21,c22,c23}.
On the other hand, some authors have also focused on the solvability
of fractional differential equations on the infinite intervals, some
excellent results were obtained, see
\cite{c24,c25,c26,c27,c28,c29,c30,c31,c32,c33,c34,c35,c36}.

 In \cite{c27} by applying standard fixed point theorems, the authors
 obtained the existence and uniqueness of solutions for
a coupled system of fractional differential equations with m-point
fractional boundary conditions
$$ \left\{
      \begin{array}{ll}
       D^{p}u(t)+f(t,v(t))=0, \ p\in (2,3),\\
       D^{q}v(t)+g(t,u(t))=0, \ q\in (2,3),\\
       u(0)=u'(0)=0,\ D^{p-1}u(+\infty)=\sum_{i=1}^{m-2}\beta_{i}
       u(\xi_{i}),\\
      v(0)=v'(0)=0,\ D^{q-1}v(+\infty)=\sum_{i=1}^{m-2}\gamma_{i} v(\xi_{i}),
     \end{array}
   \right.$$
where $t\in J=[0,+\infty),f,g\in
C(J\times\mathbb{R},\mathbb{R}),0<\xi_{1}<\xi_{2}<\cdots<\xi_{m-2}<+\infty,
\beta_{i},\gamma_{i}>0$, such that $0<\sum_{i=1}^{m-2}\beta_{i}
       u(\xi_{i})<\Gamma(p)$ and $0<\sum_{i=1}^{m-2}\gamma_{i}
       v(\xi_{i})<\Gamma(q)$,
$D^{p},D^{q}$ are the Riemann-Liouville fractional derivatives.
\par In \cite{c30} Zhai and Ren studied a coupled system of fractional differential
equations on an unbounded domain:
\begin{equation}\label{fds1} \left\{
      \begin{array}{ll}
       D^{\alpha}u(t)+\varphi(t,v(t),D^{\gamma_{1}}v(t))=0, \ \alpha\in (2,3],\gamma_{1}\in(0,1),\\
       D^{\beta}v(t)+\psi(t,u(t),D^{\gamma_{2}}u(t))=0, \ \beta\in (2,3],\gamma_{2}\in(0,1),\\
       I^{3-\alpha}u(0)=0,\
       D^{\alpha-2}u(0)=\displaystyle\int_{0}^{h}g_{1}(s)u(s)\mbox{d}s, \
       D^{\alpha-1}u(+\infty)=M
       u(\xi)+a,\\
       I^{3-\beta}v(0)=0,\
       D^{\beta-2}v(0)=\displaystyle\int_{0}^{h}g_{2}(s)v(s)\mbox{d}s, \
       D^{\beta-1}v(+\infty)=N
       v(\eta)+b,
     \end{array}
   \right.\end{equation}
where $t\in J=[0,+\infty),\varphi,\psi\in
C(J\times\mathbb{R}\times\mathbb{R}$,$J),M,N$ are real numbers
satisfying
$0<M\xi^{\alpha-1}<\Gamma(\alpha),0<N\eta^{\beta-1}<\Gamma(\beta),\xi,\eta,
h>0$, and $a,b\in \mathbb{R^{+}},g_{1},g_{2}\in L^{1}[0,h]$ are
nonnegative functions. By applying fixed point theorems, sufficient
conditions for the existence and uniqueness of solutions to the
system  \eqref{fds1} are provided , which is a natural expansion of
the results in \cite{c28}.

 In \cite{c33} Zhang et al. applied  a monotone iterative method to study
 a nonlinear fractional boundary value problem on a half line
$$ \left\{
      \begin{array}{ll}
       D^{\alpha}u(t)+f(t,u(t),D^{\alpha-1}u(t))=0, \ \alpha\in (1,2],\\
       u(0)=0,\ D^{\alpha-1}u(+\infty)=\beta u(\xi),\ \beta>0,
     \end{array}
   \right.$$
where $t\in J=[0,+\infty),f\in
C(J\times\mathbb{R}\times\mathbb{R},\mathbb{R})$. The positive
extremal solutions and iterative sequence for approximating them are
derived. A similar approach is used in \cite{c37,c38,c39,c40,c41}.
\par
Motivated by the mentioned papers, an interesting and a nature
question is if we know the existence of solution for the system
\eqref{fds}, how can we seek it? This thought motivates the research
of iterative schemes of positive solutions for the system
\eqref{fds}.
\par  By using the monotone iterative method, in this paper we
 establish two explicit monotone iterative schemes for
approximating the extreme positive solutions and construct an
explicit iterative schemes for approximating the unique positive
solution, which are more interesting and meaningful than the
traditional design route that obtains the existence of solutions.
Here we obtain not only the existence of the solution for the
system, but also the iterative schemes of the solution. Furthermore,
we extend the iterative solution problem of a single equation to the
system which is different from
\cite{c11,c26,c30,c34,c37,c38,c39,c40,c41}. Finally, the main
results extend the fractional derivative from the low-order to the
high-order fractional derivatives.

\section{Preliminaries}

We first  introduce the hypotheses that will play an important role
in subsequent proof.\\
 (H$_{1}$) \ \ $h_{i}(t)\in L[0,+\infty)$ and
$\displaystyle\int^{+\infty}_{0}h_{i}(t)t^{\alpha_{i}-1}\mbox{d}t=\Lambda _{i}<\Gamma(\alpha_{i}),f_{i}(t,0,0,0,0)\not\equiv 0,\forall t\in J,\ i=1,2.$\\
 (H$_{2}$) \ \ The nonnegative functions $a_{i0}(t),a_{ik}(t)\in L[0,+\infty)$ and
constants $\lambda_{ik}\geq 0$ satisfy
$$|f_{i}(t,u_{1},u_{2},u_{3},u_{4})|\leq a_{i0}(t)+\sum_{k=1}^{4}a_{ik}(t)|u_{k}|^{\lambda_{ik}},
\forall t\in J,\ u_{k}\in \mathbb{R}, i=1,2,\ k=1,2,3,4.$$ and
$$\displaystyle\int^{+\infty}_{0}a_{i0}(t)\mbox{d}t=a_{i0}^{\ast}<+\infty,
\displaystyle\int^{+\infty}_{0}a_{i3}(t)\mbox{d}t=a_{i3}^{\ast}<+\infty,
\displaystyle\int^{+\infty}_{0}a_{i4}(t)\mbox{d}t=a_{i4}^{\ast}<+\infty,$$
$$\int^{+\infty}_{0}a_{i1}(t)(1+t^{\alpha_{1}-1})^{\lambda_{i1}}\mbox{d}t=a_{i1}^{\ast}<+\infty,
\displaystyle\int^{+\infty}_{0}a_{i2}(t)(1+t^{\alpha_{2}-1})^{\lambda_{i2}}\mbox{d}t=a_{i2}^{\ast}<+\infty,
i=1,2.$$
 (H$_{3}$)\ \ The nonnegative functions $b_{ik}(t)\in L[0,+\infty)$ satisfy
$$
|f_{i}(t,u_{1},u_{2},u_{3},u_{4})-f_{i}(t,\bar{u}_{1},\bar{u}_{2},\bar{u}_{3},\bar{u}_{4})|\leq
\sum_{k=1}^{4}b_{ik}(t)|u_{k}-\bar{u}_{k}|,$$$$ \forall t\in
J,~u_{k},\bar{u}_{k}\in \mathbb{R},\ i=1,2,\ k=1,2,3,4.
$$
and
$$
\int^{+\infty}_{0}b_{i1}(t)(1+t^{\alpha_{1}-1})\mbox{d}t=b_{i1}^{\ast}<+\infty,
\int^{+\infty}_{0}b_{i2}(t)(1+t^{\alpha_{2}-1})\mbox{d}t=b_{i2}^{\ast}<+\infty,$$$$
\int^{+\infty}_{0}b_{i3}(t)\mbox{d}t=b_{i3}^{\ast}<+\infty,
~\int^{+\infty}_{0}b_{i4}(t)\mbox{d}t=b_{i4}^{\ast}<+\infty,
\int^{+\infty}_{0}|f_{i}(t,0,0,0,0)|\mbox{d}t=\tau_{i}<+\infty,\
i=1,2.
$$
(H$_{4}$)\  Functions $f_{i}(t,u_{1},u_{2},u_{3},u_{4})$ are
increasing with respect to the variables
$u_{1},u_{2},u_{3},u_{4},\forall t\in J, i=1,2$.

Next we list some definitions and lemmas that are helpful to the
proof of principal theorems.

 {\bf Definition 2.1}(see \cite{bk1,bk3}).\ The
Riemann-Liouville fractional integral of order $q>0$ for an
{\color{red}integrable} function $g$ is defined as
 $$
  I^{q}g(x)=\frac{1}{\Gamma(q)}\int^{x}_{0}(x-t)^{q-1}g(t)\mathrm{d}
  t,$$
  provided that the integral exists.

 {\bf Definition 2.2}(see \cite{bk1,bk3}).\ The Riemann-Liouville fractional derivative of
   order $q>0$ for an {\color{red}integrable} function $g$ is defined as
   $$
  D^{q} g(x)=\frac{1}{\Gamma(n-q)}
  \Big(\frac{d}{dx}\Big)^{n}\int^{x}_{0}(x-t)^{n-q-1}g(t)\mathrm{d} t,
   $$
where  $n=[q]+1, [\alpha]$ is the smallest integer greater than or
equal to $\alpha$, provided that the right-hand side is pointwise
defined on $(0,+\infty)$.

  {\bf Lemma 2.1}(see \cite{bk1,bk3}). \ Let $q>0$ and $u\in
C(0,1)\cap L(0, 1)$. Then the general solution of fractional
differential equation $D^{q}u(t)=0$ is
$$ u(t)=c_{1}t^{q-1}+c_{2}t^{q-2}+\cdots+c_{n}t^{q-n},
$$
where $c_{i}\in \mathbb{R},i=1,2,\cdots,n$ and $n-1<q<n$.

 {\bf Lemma 2.2.} \  Let $y_{i}\in C[0,+\infty)$ with
$\displaystyle\int^{+\infty}_{0}h_{i}(t)t^{\alpha_{i}-1}\mbox{d}t\neq\Gamma(\alpha_{i}),n_{i}-1<\alpha_{i}\leq
n_{i},i=1,2$. Then the fractional differential system boundary value
problem
\begin{equation}\label{fdsbvp} \left\{
      \begin{array}{ll}
       D^{\alpha_{1}}u(t)+y_{1}(t)=0,\
       n_{1}-1<\alpha_{1}\leq n_{1},\\
        D^{\alpha_{2}}v(t)+y_{2}(t)=0,\
        n_{2}-1<\alpha_{2}\leq n_{2},\\
       u(0)=u'(0)=\cdots=u^{(n_{1}-2)}(0)=0,D^{\alpha_{1}-1}u(+\infty)=
       \displaystyle\int_{0}^{+\infty}h_{1}(t)u(t)\mbox{d}t,\\
        v(0)=v'(0)=\cdots=v^{(n_{2}-2)}(0)=0,D^{\alpha_{2}-1}v(+\infty)=
       \displaystyle\int_{0}^{+\infty}h_{2}(t)v(t)\mbox{d}t,\\
     \end{array}
   \right.\end{equation}
 has the integral representation
\begin{equation}\label{fdsbvp1}\left\{
      \begin{array}{ll}
    u(t)= \displaystyle\int^{+\infty}_{0}K_{1}(t,s)y_{1}(s)\mbox{d}s,\\
     v(t)= \displaystyle\int^{+\infty}_{0}K_{2}(t,s)y_{2}(s)\mbox{d}s,
 \end{array}
   \right.\end{equation}
where
\begin{equation}\label{k}K_{i}(t,s)=K_{i1}(t,s)+K_{i2}(t,s),i=1,2.\end{equation}
with
\begin{equation}\label{k1}K_{i1}(t,s)=\frac{1}{\Gamma(\alpha_{i})}\left\{
      \begin{array}{ll}
      t^{\alpha_{i}-1}-(t-s)^{\alpha_{i}-1},0\leq s\leq t\leq +\infty, \\
    t^{\alpha_{i}-1},0\leq t\leq s\leq +\infty,\\
     \end{array}
   \right.\end{equation}
   \begin{equation}\label{k2}K_{i2}(t,s)=
       \frac{t^{\alpha_{i}-1}}{\Gamma(\alpha_{i})-\Lambda_{i}}\int^{+\infty}_{0}h_{i}(t)K_{i1}(t,s)\mbox{d}t.
    \end{equation}

{\bf Proof.} \ From Lemma 2.1, we can turn differential system
\eqref{fdsbvp} into an equivalent integral system
\begin{equation}\label{eis}\left\{
      \begin{array}{ll}
    u(t)=-I^{\alpha_{1}}y_{1}(t)+ c_{11}t^{\alpha_{1}-1}+c_{12}t^{\alpha_{1}-2}+\ldots+c_{1n_{1}}t^{\alpha_{1}-n_{1}},\\
     v(t)=-I^{\alpha_{2}}y_{2}(t)+
     c_{21}t^{\alpha_{2}-1}+c_{22}t^{\alpha_{2}-2}+\ldots+c_{2n_{2}}t^{\alpha_{2}-n_{2}},
 \end{array}
   \right. \end{equation}
where
$c_{11},c_{12},\cdots,c_{1n_{1}},c_{21},c_{22},\cdots,c_{2n_{2}}$
are arbitrary constants.
 With the help of conditions $u(0)=u'(0)=\cdots=u^{(n_{1}-2)}(0)=0$ and $v(0)=v'(0)=\cdots=v^{(n_{2}-2)}(0)=0$,
it is easy to know that
$c_{12}=c_{13}=\cdots=c_{1n_{1}}=c_{22}=c_{23}=\cdots=c_{2n_{2}}=0$.
From  \eqref{eis} we have
\begin{equation}\label{eis1}\left\{
      \begin{array}{ll}
    u(t)=-\displaystyle\frac{1}{\Gamma(\alpha_{1})}\int_{0}^{t}(t-s)^{\alpha_{1}}y_{1}(s)\mbox{d}s+ c_{11}t^{\alpha_{1}-1},\\
     v(t)=-\displaystyle\frac{1}{\Gamma(\alpha_{2})}\int_{0}^{t}(t-s)^{\alpha_{2}}y_{2}(s)\mbox{d}s+ c_{21}t^{\alpha_{2}-1}.
 \end{array}
   \right.\end{equation}
Then
\begin{equation}\label{eis2}\left\{
      \begin{array}{ll}
     D^{\alpha_{1}-1}u(t)= c_{11}\Gamma(\alpha_{1})-\displaystyle\int_{0}^{t}y_{1}(s)\mbox{d}s,\\
     D^{\alpha_{2}-1}v(t)=c_{21}\Gamma(\alpha_{2})-\displaystyle\int_{0}^{t}y_{2}(s)\mbox{d}s.
 \end{array}
   \right.\end{equation}
Hence
 \begin{equation}\label{eis3}\left\{
      \begin{array}{ll}
     D^{\alpha_{1}-1}u(+\infty)= c_{11}\Gamma(\alpha_{1})-\displaystyle\int_{0}^{+\infty}y_{1}(s)\mbox{d}s,\\
     D^{\alpha_{2}-1}v(+\infty)=c_{21}\Gamma(\alpha_{2})-\displaystyle\int_{0}^{+\infty}y_{2}(s)\mbox{d}s.
 \end{array}
   \right.\end{equation}
 Based on the conditions $D^{\alpha_{1}-1}u(+\infty)=
       \int_{0}^{+\infty}h_{1}(t)u(t)\mbox{d}t$ and $D^{\alpha_{2}-1}v(+\infty)=
       \int_{0}^{+\infty}h_{2}(t)v(t)\mbox{d}t$, we have
 \begin{equation}\label{c1}\left\{
      \begin{array}{ll}
      c_{11}=\displaystyle\frac{1}{\Gamma(\alpha_{1})}\int_{0}^{+\infty}h_{1}(t)u(t)\mbox{d}t
      +\frac{1}{\Gamma(\alpha_{1})}\int_{0}^{+\infty}y_{1}(s)\mbox{d}s,\\
     c_{21}=\displaystyle\frac{1}{\Gamma(\alpha_{2})}\int_{0}^{+\infty}h_{2}(t)v(t)\mbox{d}t
     +\frac{1}{\Gamma(\alpha_{2})}\int_{0}^{+\infty}y_{2}(s)\mbox{d}s.
 \end{array}
   \right.\end{equation}
 Submitting \eqref{c1} to \eqref{eis2}, we know
\begin{equation}\label{c2}\left\{
      \begin{array}{ll}
      u(t)&=-\displaystyle\frac{1}{\Gamma(\alpha_{1})}\int_{0}^{t}(t-s)^{\alpha_{1}}y_{1}(s)\mbox{d}s+
      \frac{t^{\alpha_{1}-1}}{\Gamma(\alpha_{1})}\Big[\int_{0}^{+\infty}h_{1}(t)u(t)\mbox{d}t+\int_{0}^{+\infty}y_{1}(s)\mbox{d}s\Big]\\
       & =\displaystyle\int_{0}^{\infty}K_{11}(t,s)y_{1}(s)\mbox{d}s+\frac{t^{\alpha_{1}-1}}{\Gamma(\alpha_{1})}\int_{0}^{+\infty}h_{1}(t)u(t)\mbox{d}t,\\
     v(t)&=-\displaystyle\frac{1}{\Gamma(\alpha_{2})}\int_{0}^{t}(t-s)^{\alpha_{2}}y_{2}(s)\mbox{d}s+
      \frac{t^{\alpha_{2}-1}}{\Gamma(\alpha_{2})}\Big[\int_{0}^{+\infty}h_{2}(t)v(t)\mbox{d}t+\int_{0}^{+\infty}y_{2}(s)\mbox{d}s\Big]\\
       & =\displaystyle\int_{0}^{\infty}K_{21}(t,s)y_{2}(s)\mbox{d}s+\frac{t^{\alpha_{2}-1}}{\Gamma(\alpha_{2})}\int_{0}^{+\infty}h_{2}(t)v(t)\mbox{d}t.\\
 \end{array}
   \right.\end{equation}
Multiplying both sides of the above equality  by $h_{1}(t)$ and
$h_{2}(t)$ and integrating from 0 to $+\infty$ , we obtain
$$\left\{
      \begin{array}{ll}
     \displaystyle\int_{0}^{+\infty}h_{1}(t)u(t)\mbox{d}t=
     \frac{\Gamma(\alpha_{1})}{\Gamma(\alpha_{1})-\Lambda_{1}}\int_{0}^{+\infty}h_{1}(t)\int_{0}^{+\infty}K_{11}(t,s)y_{1}(s)\mbox{d}s\mbox{d}t,\\
    \displaystyle\int_{0}^{+\infty}h_{2}(t)v(t)\mbox{d}t =
    \frac{\Gamma(\alpha_{2})}{\Gamma(\alpha_{2})-\Lambda_{2}}\int_{0}^{+\infty}h_{2}(t)\int_{0}^{+\infty}K_{21}(t,s)y_{2}(s)\mbox{d}s\mbox{d}t.
 \end{array}
   \right.$$
 Combining \eqref{c2}, we have
$$\left\{
      \begin{array}{ll}
      u(t)&=\displaystyle\int_{0}^{+\infty}K_{11}(t,s)y_{1}(s)\mbox{d}s
      +\frac{t^{\alpha_{1}-1}}{\Gamma(\alpha_{1})-\Lambda_{1}}\int_{0}^{+\infty}h_{1}(t)\int_{0}^{+\infty}K_{11}(t,s)y_{1}(s)\mbox{d}s\mbox{d}t\\
      &=\displaystyle\int_{0}^{\infty}K_{11}(t,s)y_{1}(s)\mbox{d}s
      +\int_{0}^{\infty}K_{12}(t,s)y_{1}(s)\mbox{d}s,\\
      &=\displaystyle\int_{0}^{\infty}K_{1}(t,s)y_{1}(s)\mbox{d}s,\\
     v(t)&=\displaystyle\int_{0}^{+\infty}K_{21}(t,s)y_{2}(s)\mbox{d}s+
     \frac{t^{\alpha_{2}-1}}{\Gamma(\alpha_{2})-\Lambda_{2}}\int_{0}^{+\infty}h_{2}(t)\int_{0}^{+\infty}K_{21}(t,s)y_{2}(s)\mbox{d}s\mbox{d}t\\
     &=\displaystyle\int_{0}^{+\infty}K_{21}(t,s)y_{2}(s)\mbox{d}s+
    \int_{0}^{+\infty}K_{22}(t,s)y_{2}(s)\mbox{d}s\\
    &=\displaystyle\int_{0}^{+\infty}K_{2}(t,s)y_{2}(s)\mbox{d}s.
 \end{array}
   \right.$$
The proof is completed.

 {\bf Remark 2.1.} \ From \eqref{fdsbvp1}, \eqref{k}, \eqref{k1} and
\eqref{k2}, by direct calculation, we have
$$\left\{
      \begin{array}{ll}
    D^{\alpha_{1}-1}u(t)=\displaystyle\int_{0}^{+\infty}K_{1}^{\ast}(t,s)y_{1}(s)\mbox{d}s,\\
     D^{\alpha_{2}-1}v(t)= \displaystyle\int_{0}^{+\infty}K_{2}^{\ast}(t,s)y_{2}(s)\mbox{d}s,
 \end{array}
   \right.$$
 where
$$K^{\ast}_{i}(t,s)=K^{\ast}_{i1}(t,s)+K^{\ast}_{i2}(t,s),\ i=1,2.$$ with
$$ K^{\ast}_{i1}(t,s)=\left\{
      \begin{array}{ll}
      0,0\leq s\leq t\leq +\infty, \\
    1,0\leq t\leq s\leq +\infty,
     \end{array}\right.\\
    K^{\ast}_{i2}(t,s)=
       \frac{\Gamma(\alpha_{i})}{\Gamma(\alpha_{i})-\Delta_{i}}\int^{+\infty}_{0}h_{i}(t)K_{i1}(t,s)\mbox{d}t.$$

 {\bf Lemma 2.3.} \ For $(s,t)\in J\times J$, if hypothesis (H$_{1}$) is satisfied, then
  $$0\leq
  K_{i}(t,s)\leq\frac{t^{\alpha_{i}-1}}{\Gamma(\alpha_{i})-\Lambda_{i}},\ \
 0\leq
 \frac{K_{i}(t,s)}{1+t^{\alpha_{i}-1}}\leq\frac{1}{\Gamma(\alpha_{i})-\Lambda_{i}},\ i=1,2.$$

 {\bf Proof.} \ From \eqref{k1} and \eqref{k2}, it is obvious that
$$0\leq K_{i1}(t,s)\leq
 \frac{t^{\alpha_{i}-1}}{\Gamma(\alpha_{i})},\forall
(t,s)\in J\times
 J,$$ and $$0\leq K_{i2}(t,s)\leq
 \frac{t^{\alpha_{i}-1}}{\Gamma(\alpha_{i})-\Lambda_{i}}\int^{+\infty}_{0} \frac{h_{i}(t)t^{\alpha_{i}-1}}{\Gamma(\alpha_{i})}\mbox{d}t
 =\frac{\Lambda_{i} t^{\alpha_{i}-1}}{\Gamma(\alpha_{i})(\Gamma(\alpha_{i})-\Lambda_{i})},\forall
(t,s)\in J\times
 J.$$
 So
 $$0\leq K_{i}(t,s)=K_{i1}(t,s)+K_{i2}(t,s)\leq \frac{ t^{\alpha_{i}-1}}{\Gamma(\alpha_{i})-\Lambda_{i}},\forall
(t,s)\in J\times
 J.$$
Furthermore
 $$0\leq
 \frac{K_{i}(t,s)}{1+t^{\alpha_{i}-1}}\leq\frac{1}{\Gamma(\alpha_{i})-\Lambda_{i}},\forall
(t,s)\in J\times
 J.$$
The proof is completed.

 {\bf Remark 2.2.} \ From Remark 2.1, by direct calculation, we can easily
 know that
 $$0\leq K_{i}^{\ast}(t,s)=K^{\ast}_{i1}(t,s)+K^{\ast}_{i2}(t,s)\leq 1+\frac{ \Lambda_{i}}{\Gamma(\alpha)-\Lambda_{i}}
 =\frac{\Gamma(\alpha_{i})}{\Gamma(\alpha_{i})-\Lambda_{i}},
 \forall
(t,s)\in J\times
 J,\ i=1,2,$$

 Let $E=\{u\in C(J,\mathbb{R})|\sup_{t\in
 J}\frac{|u(t)|}{1+t^{\alpha_{1}-1}}<+\infty\}$ and $X=\{u\in E,D^{\alpha_{1}-1}u\in C(J,\mathbb{R})| \ \sup_{t\in
 J}|D^{\alpha_{1}-1}u(t)|<+\infty\}$ be equipped with the norm
 $$\|u\|_{X}=\max\{\|u\|_{0},
\|D^{\alpha_{1}-1}u\|_{1}\},$$ where $\|u\|_{0}=\sup_{t\in
 J}\frac{|u(t)|}{1+t^{\alpha_{1}-1}}$  and  $\|D^{\alpha_{1}-1}u\|_{1}=\sup_{t\in
 J}|D^{\alpha_{1}-1}u(t)|$. Also let
$F=\{v\in C(J,\mathbb{R})|\sup_{t\in
 J}\frac{|v(t)|}{1+t^{\alpha_{2}-1}}<+\infty\}$ and $Y=\{v\in F,D^{\alpha_{2}-1}v\in C(J,\mathbb{R})|\sup_{t\in
 J}|D^{\alpha_{2}-1}v(t)|<+\infty\}$ be equipped with the norm
 $$\|u\|_{Y}=\max\{\|v\|_{0},
\|D^{\alpha_{2}-1}v\|_{1}\},$$ where $\|v\|_{0}=\sup_{t\in
 J}\frac{|v(t)|}{1+t^{\alpha_{2}-1}}$  and  $\|D^{\alpha_{2}-1}v\|_{1}=\sup_{t\in
 J}|D^{\alpha_{2}-1}v(t)|$. Thus the space $(X,\|\cdot\|_{X})$ and
 $(Y,\|\cdot\|_{Y})$ are two Banach spaces which have been shown in
 \cite{c24}. Moreover, the product space $(X\times Y,\|\cdot\|_{X\times Y})$ is also a
 Banach space with the norm $$\|\cdot\|_{X\times Y}=\max\{\|u\|_{X},
\|v\|_{Y}\}.$$

{\bf  Lemma 2.4.} \ If hypothesis (H$_{2}$) is satisfied, then for
$\forall (u,v)\in X\times Y$, we have
$$\int_{0}^{+\infty}|f_{i}(s,u(s),v(s),D^{\alpha_{1}-1}u(s),D^{\alpha_{2}-1}v(s))|\mbox{d}s
\leq a_{i0}^{\ast}+\sum_{k=1}^{4}a_{ik}^{\ast}||(u,v)||_{X\times
Y}^{\lambda_{ik}},i=1,2.$$

 {\bf Proof.} For $\forall (u,v)\in
X\times Y$, by hypothesis (H$_{2}$), we have
\[\aligned &\int_{0}^{+\infty}|f_{i}(s,u(s),v(s),D^{\alpha_{1}-1}u(s),D^{\alpha_{2}-1}v(s))|\mbox{d}s\\
\leq&\int_{0}^{+\infty}
 \Big(a_{i0}(s)+a_{i1}(s)|u(s)|^{\lambda_{i1}}+a_{i2}(s)|v(s)|^{\lambda_{i2}}+a_{i3}(s)|D^{\alpha_{1}-1}u(s))|^{\lambda_{i3}}+a_{i4}(s)|D^{\alpha_{2}-1}v(s))|^{\lambda_{i4}}\Big)\mbox{d}s
 \\
 \leq&a_{i0}^{\ast}+\int_{0}^{+\infty}a_{i1}(s))(1+s^{\alpha_{1}-1})^{\lambda_{i1}}\frac{|u(s)|^{\lambda_{i1}}}
 {(1+s^{\alpha_{1}-1})^{\lambda_{i1}}}ds+\int_{0}^{+\infty}a_{i2}(s))(1+s^{\alpha_{2}-1})^{\lambda_{i2}}\frac{|v(s)|^{\lambda_{i2}}}{(1+s^{\alpha_{2}-1})^{\lambda_{i2}}}ds\\
 &+\int_{0}^{+\infty}a_{i3}(s)|D^{\alpha_{1}-1}u(s)|^{\lambda_{i3}}ds
 +\int_{0}^{+\infty}a_{i4}(s)|D^{\alpha_{2}-1}v(s)|^{\lambda_{i4}}ds\\
\leq&
a_{i0}^{\ast}+a_{i1}^{\ast}||u||_{X}^{\lambda_{i1}}+a_{i2}^{\ast}||v||_{Y}^{\lambda_{i2}}
+a_{i3}^{\ast}||u||_{X}^{\lambda_{i3}}+a_{i4}^{\ast}||v||_{Y}^{\lambda_{i4}}\\
\leq& a_{i0}^{\ast}+\sum_{k=1}^{4}a_{ik}^{\ast}||(u,v)||_{X\times
Y}^{\lambda_{ik}},\ i=1,2.
\endaligned\]

{\bf Lemma 2.5.} \ If hypothesis (H$_{3}$) is satisfied, then for
$\forall (u,v)\in X\times Y$, we have
$$\int_{0}^{+\infty}|f_{i}(s,u(s),v(s),D^{\alpha_{1}-1}u(s),D^{\alpha_{2}-1}v(s))|\mbox{d}s
\leq \sum_{k=1}^{4}b_{ik}^{\ast}||(u,v)||_{X\times Y}+\tau_{i},\
i=1,2.$$

 { \bf Proof. } For $\forall (u,v)\in X\times Y$, by
hypothesis (H$_{3}$), we have
\[\aligned &\int_{0}^{+\infty}|f_{i}(s,u(s),v(s),D^{\alpha_{1}-1}u(s),D^{\alpha_{2}-1}v(s))|\mbox{d}s\\
 =&\int_{0}^{+\infty}|f_{i}(s,u(s),v(s),D^{\alpha_{1}-1}u(s),D^{\alpha_{2}-1}v(s))-f_{i}(s,0,0,0,0)+f_{i}(s,0,0,0,0)|\mbox{d}s\\
 \leq&\int_{0}^{+\infty}|f_{i}(s,u(s),v(s),D^{\alpha_{1}-1}u(s),D^{\alpha_{2}-1}v(s))
 -f_{i}(s,0,0,0,0)|\mbox{d}s+\int_{0}^{+\infty}|f_{i}(s,0,0,0,0)|\mbox{d}s\\
\leq&
\int_{0}^{+\infty}b_{i1}(s)(1+s^{\alpha_{1}-1})\frac{|u(s)|}{1+s^{\alpha_{1}-1}}\mbox{d}s+
\int_{0}^{+\infty}b_{i2}(s)(1+s^{\alpha_{2}-1})\frac{|v(s)|}{1+s^{\alpha_{2}-1}}\mbox{d}s\\
&+\int_{0}^{+\infty}b_{i3}(s)|D^{\alpha_{1}-1}u(s)|\mbox{d}s
+\int_{0}^{+\infty}b_{i4}(s)|D^{\alpha_{2}-1}v(s)|\mbox{d}s
+\int_{0}^{+\infty}|f_{i}(s,0,0,0,0)|\mbox{d}s\\
\leq&
b_{i1}^{\ast}||u||_{X}+b_{i2}^{\ast}||v||_{Y}+b_{i3}^{\ast}||u||_{X}+b_{i4}^{\ast}||v||_{Y}
+\tau_{i}\\
\leq& \sum_{k=1}^{4}b_{ik}^{\ast}||(u,v)||_{X\times Y}+\tau_{i},\
i=1,2.\endaligned\]

{\bf Lemma 2.6} (see \cite{c24}). \ Let $U\subset X$ be a bounded
set. Then $U$ is a
relatively compact in $X$ if the following conditions hold:\\
 (i)\ For any $u\in U,\displaystyle\frac{u(t)}{1+t^{\alpha-1}}$ and $D^{\alpha-1}u(t)$
are equicontinuous on any compact
interval of J;\\
 (ii)\ For any $\varepsilon>0$, there is a
constant $C=C(\varepsilon)>0$ such that
$\displaystyle|\frac{u(t_{1})}{1+t_{1}^{\alpha-1}}-\frac{u(t_{2})}{1+t_{2}^{\alpha-1}}|<\varepsilon$
and $|D^{\alpha-1}u(t_{1})-D^{\alpha-1}u(t_{2})|<\varepsilon$ for
any $t_{1},t_{2}\geq C$ and $u\in U$.

 We define the cone $P\subset X\times Y$ by $P=\{(u,v)\in
X\times Y| u(t)\geq0,v(t)\geq0,
 D^{\alpha_{1}-1}u(t)\geq0,D^{\alpha_{2}-1}v(t)\geq0,t\in J\}$. By
 Lemma 2.2, let $T:P\rightarrow P$ be the operator defined as

 \begin{equation}\label{op}
  T(u,v)(t)=
\left(\aligned
&{T_{1}(u,v)(t) } \\
&{  T_{2}(u,v)(t)}
\endaligned\right)=
\left(\aligned
&{\int^{+\infty}_{0}K_{1}(t,s)f_{1}(s,u(s),v(s),D^{\alpha_{1}-1}u(s),D^{\alpha_{2}-1}v(s))\mbox{d}s } \\
&{\int^{+\infty}_{0}K_{2}(t,s)f_{2}(s,u(s),v(s),D^{\alpha_{1}-1}u(s),D^{\alpha_{2}-1}v(s))\mbox{d}s
}
\endaligned\right)
   \end{equation}
 By Remark 2.1, we also define
 \begin{equation}\label{op1}
\left(\aligned
&{ D^{\alpha_{1}-1} T_{1}(u,v)(t) } \\
&{  D^{\alpha_{2}-1}T_{2}(u,v)(t)}
\endaligned\right)=
\left(\aligned
&{\int_{0}^{+\infty}K_{1}^{\ast}(t,s)f_{1}(s,u(s),v(s),
 D^{\alpha_{1}-1}u(s),D^{\alpha_{2}-1}v(s))\mbox{d}s  } \\
&{\int_{0}^{+\infty}K_{2}^{\ast}(t,s)f_{2}(s,u(s),v(s),D^{\alpha_{1}-1}u(s),D^{\alpha_{2}-1}v(s))\mbox{d}s
}
\endaligned\right)
 \end{equation}
It is easy to know that the system \eqref{fds} has a solution
 if and only if the operator equation $(u,v)=T(u,v)$ has a fixed point,
 where $T$ is given by \eqref{op}. In fact, if $(u,v)$ is a solution for the system \eqref{fds}, by lemma 2.2, we can obtain
  $$
      \begin{array}{ll}
      u=\displaystyle\int^{+\infty}_{0}K_{1}(t,s)f_{1}(s,u(s),v(s),D^{\alpha_{1}-1}u(s),D^{\alpha_{2}-1}v(s))\mbox{d}s=T_{1}(u,v),\\
     v=\displaystyle\int^{+\infty}_{0}K_{2}(t,s)f_{2}(s,u(s),v(s),D^{\alpha_{1}-1}u(s),D^{\alpha_{2}-1}v(s))\mbox{d}s=T_{2}(u,v).
 \end{array}
 $$
 That is, $(u,v)$ is a fixed point for the operator equation $(u,v)=T(u,v)$. On the contrary, the
Riemann-Liouville fractional
 derivation on both sides of the operator equation is
 $$
      \begin{array}{ll}
     D^{\alpha_{1}}u(t)=D^{\alpha_{1}} T_{1}(u,v)(t)=-f_{1}(s,u(s),v(s),D^{\alpha_{1}-1}u(s),D^{\alpha_{2}-1}v(s)),\\
     D^{\alpha_{2}}v(t)=D^{\alpha_{2}}T_{2}(u,v)(t)=-f_{2}(s,u(s),v(s),D^{\alpha_{1}-1}u(s),D^{\alpha_{2}-1}v(s)).
 \end{array}
 $$
 Combining \eqref{eis3} and \eqref{c1}, we can obtain
$$D^{\alpha_{1}-1}u(+\infty)=
       \displaystyle\int_{0}^{+\infty}h_{1}(t)u(t)\mbox{d}t,\ \
        D^{\alpha_{2}-1}v(+\infty)=
       \displaystyle\int_{0}^{+\infty}h_{2}(t)v(t)\mbox{d}t,$$
That is, $(u,v)$ is a solution for the system \eqref{fds}.

 {\bf Lemma 2.7.}\ If the hypotheses (H$_{1})$ and (H$_{2})$ are satisfied, then the
operator $T:P\rightarrow P$ is completely
continuous.\\
{ \bf Proof.} \ First it is easy to know $T:P\rightarrow P$. Since
$K_{i}(t,s)\geq 0$ and $f_{i}\geq 0$, we have $T_{i}(u,v)(t)\geq
0,\forall (u,v)\in P,\ t\in J,\ i=1,2$.

Next we prove in three steps that the operator $T:P\rightarrow P$ is
relatively compact.

\textbf{Step 1} Let $U=\{(u,v)|(u,v)\in P,||(u,v)||_{X\times Y}\leq
M\}$. For $\forall (u,v)\in U$, by Lemma 2.3, Remark 2.2 and Lemma
2.4,
    we obtain
\begin{equation}\label{y0}\begin{array}{rcl}
 ||T_{1}(u,v)||_{0}&=&\displaystyle\sup_{t\in J}\Big|\int_{0}^{+\infty}\frac{K_{1}(t,s)}{1+t^{\alpha_{1}-1}}
 f_{1}(s,u(s),v(s),D^{\alpha_{1}-1}u(s),D^{\alpha_{2}-1}v(s))\mbox{d}s\Big|\\
&\leq&\displaystyle\frac{1}{\Gamma(\alpha_{1})-\Lambda_{1}}\int_{0}^{+\infty}
 |f_{1}(s,u(s),v(s),D^{\alpha_{1}-1}u(s),D^{\alpha_{2}-1}v(s))|\mbox{d}s
 \\ &\leq&\displaystyle\frac{1}{\Gamma(\alpha_{1})-\Lambda_{1}}\Big[a_{10}^{\ast}+\sum_{k=1}^{4}a_{1k}^{\ast}||(u,v)||_{X\times
Y}^{\lambda_{1k}}\Big]
\end{array}\end{equation}
and
\begin{equation}\label{y00}\begin{array}{rcl}
\displaystyle ||T_{1}(u,v)||_{1}&=&\displaystyle\sup_{t\in
J}\Big|\int_{0}^{\infty}K_{1}^{\ast}(t,s)f_{1}(s,u(s),v(s),D^{\alpha_{1}-1}u(s),D^{\alpha_{2}-1}v(s))\mbox{d}s\Big|\\
 &\leq&\displaystyle\frac{\Gamma(\alpha_{1})}{\Gamma(\alpha_{1})-\Lambda_{1}}\int_{0}^{+\infty}|f_{1}(s,u(s),v(s),D^{\alpha_{1}-1}u(s),D^{\alpha_{2}-1}v(s))|\mbox{d}s\\
 &\leq&\displaystyle\frac{\Gamma(\alpha_{1})}{\Gamma(\alpha_{1})-\Lambda_{1}}\Big[a_{10}^{\ast}+\sum_{k=1}^{4}a_{1k}^{\ast}||(u,v)||_{X\times
Y}^{\lambda_{1k}}\Big].
\end{array}\end{equation}
Thus
$$||T_{1}(u,v)||_{X}=\max\Big\{\|T_{1}(u,v)\|_{0},\|T_{1}(u,v)\|_{1}\Big\}\leq
\frac{\max\{1,\Gamma(\alpha_{1})\}}{\Gamma(\alpha_{1})-\Lambda_{1}}
\Big[a_{10}^{\ast}+\sum_{k=1}^{4}a_{1k}^{\ast}M^{\lambda_{1k}}\Big].$$
Similarly
$$||T_{2}(u,v)||_{Y}=\max\Big\{\|T_{2}(u,v)\|_{0},\|T_{2}(u,v)\|_{1}\Big\}\leq
\frac{\max\{1,\Gamma(\alpha_{2})\}}{\Gamma(\alpha_{2})-\Lambda_{2}}
\Big[a_{20}^{\ast}+\sum_{k=1}^{4}a_{2k}^{\ast}M^{\lambda_{2k}}\Big].$$
Then
$$\begin{array}{rcl}
\displaystyle ||T(u,v)||_{X\times
Y}&=&\max\Big\{\|T_{1}(u,v)\|_{X},\|T_{2}(u,v)\|_{Y}\Big\}\\
 &\leq& \max\displaystyle\Big\{\frac{\max\{1,\Gamma(\alpha_{1})\}}{\Gamma(\alpha_{1})-\Lambda_{1}}
\Big(a_{10}^{\ast}+\sum_{k=1}^{4}a_{1k}^{\ast}M^{\lambda_{1k}}\Big),\\
&&\displaystyle\frac{\max\{1,\Gamma(\alpha_{2})\}}{\Gamma(\alpha_{2})-\Lambda_{2}}
\Big(a_{20}^{\ast}+\sum_{k=1}^{4}a_{2k}^{\ast}M^{\lambda_{2k}}\Big)\Big\}.
\end{array}$$
which means that $TU$ is uniformly bounded.

 \textbf{Step 2} \ Let $I\subset J$ be any compact interval. Then, for all $t_{1},t_{2}\in I,t_{2}>t_{1}$
and $(u,v)\in U$, we have
  \begin{equation}\label{y}
  \begin{array}{rcl}&&\displaystyle\Big|\frac{T_{1}(u,v)(t_{2})}{1+t^{\alpha_{1}-1}_{2}}
 -\frac{T_{1}(u,v)(t_{1})}{1+t^{\alpha_{1}-1}_{1}}\Big|\\
&\leq&\Big|\displaystyle\int_{0}^{+\infty}\Big(\frac{K_{1}(t_{2},s)}{1+t^{\alpha_{1}-1}_{2}}
-\frac{K_{1}(t_{1},s)}{1+t^{\alpha_{1}-1}_{1}}\Big) f_{1}(s,u(s),v(s),D^{\alpha_{1}-1}u(s),D^{\alpha_{2}-1}v(s))\mbox{d}s\Big|\\
&\leq&
\displaystyle\int_{0}^{+\infty}\Big|\frac{K_{1}(t_{2},s)}{1+t^{\alpha_{1}-1}_{2}}
-\frac{K_{1}(t_{1},s)}{1+t^{\alpha_{1}-1}_{1}}\Big|\big|
f_{1}(s,u(s),v(s),D^{\alpha_{1}-1}u(s),D^{\alpha_{2}-1}v(s))\big|\mbox{d}s
  \end{array}\end{equation}
Noticing that $K_{1}(t,s)/(1+t^{\alpha_{1}-1})$ is uniformly
continuous for any $(t,s)\in I\times I$. In the meantime, the
function $K_{1}(t,s)/(1+t^{\alpha_{1}-1})$ only relys on $t$ for
$s\geq t$, which infers that $K_{1}(t,s)/(1+t^{\alpha_{1}-1})$ is
uniformly continuous on $I\times (J\setminus I)$. Therefore, for all
$s\in J$ and $t_{1},t_{2}\in I$, we have
\begin{equation}\label{y1} \forall \epsilon>0, \exists \delta(\epsilon)~ such~ that~ if
~|t_{1}-t_{2}|<\delta,~ then
~\Big|\frac{K_{1}(t_{2},s)}{1+t^{\alpha_{1}-1}_{2}}
-\frac{K_{1}(t_{1},s)}{1+t^{\alpha_{1}-1}_{1}}\Big|<\epsilon.
\end{equation}
\par By Lemma 2.4, for all $(u,v)\in U$, we can obtain
$$\int_{0}^{+\infty}|f_{1}(s,u(s),v(s),D^{\alpha_{1}-1}u(s),D^{\alpha_{2}-1}v(s))|ds\leq
\Big[a_{10}^{\ast}+\sum_{k=1}^{4}a_{1k}^{\ast}M^{\lambda_{1k}}\Big]<\infty,$$
together \eqref{y} and \eqref{y1}, which means that
$T_{1}(u,v)(t)/(1+t^{\alpha_{1}-1})$ is equicontinuous on $I$.
\par Note that
 $$D^{\alpha_{1}-1}T_{1}(u,v)(t)=\int_{0}^{+\infty}K_{1}^{\ast}(t,s)f_{1}(s,u(s),v(s),D^{\alpha_{1}-1}u(s),D^{\alpha_{2}-1}v(s))\mbox{d}s$$
 and the function $K_{1}^{\ast}(t,s)\in C(J\times J)$ doesn't rely on
 $t$, which means that $D^{\alpha_{1}-1}T_{1}(u,v)(t)$ is
 equicontinuous on $I$.
In the same way, we can show that
$T_{2}(u,v)(t)/(1+t^{\alpha_{2}-1})$  and
$D^{\alpha_{2}-1}T_{2}(u,v)(t)$ are equicontinuous. Thus $T_{1}$ and
$T_{2}$ is equicontinuous on $I$.

As a natural result, the operator $T$ is equicontinuous for all
$(u,v)\in U$  on any compact interval $I$ of $J$.

\textbf{Step 3} \ We show the operator $T$ is equiconvergent at
$+\infty$.
 Since
 $$\lim_{t\rightarrow +\infty}\frac{K_{i}(t,s)}{1+t^{\alpha_{i}-1}}
 =\frac{1}{\Gamma(\alpha_{i})}+\frac{1}{\Gamma(\alpha_{i})
 -\Lambda_{i}}\int^{+\infty}_{0}h(t)K_{i1}(t,s)\mbox{d}t\leq\frac{1}{\Gamma(\alpha_{i})-\Lambda_{i}}<+\infty,\ i=1,2,$$
 by knowledge of limit theory, we can deduce that for any
 $\epsilon>0$, there exists a constant $C=C(\epsilon)>0$, for any
$t_{1},t_{2}\geq C$ and $s\in J$, such that
$$\Big|\frac{K_{i}(t_{2},s)}{1+t^{\alpha_{i}-1}_{2}}
-\frac{K_{i}(t_{1},s)}{1+t^{\alpha_{i}-1}_{1}}\Big|<\epsilon,\
i=1,2,$$ Therefore, by Lemma 2.4
 and \eqref{y}, we conclude that $T_{i}(u,v)(t)/1+t^{\alpha_{i}-1}(i=1,2)$ are equiconvergent at
 $+\infty$. As the function $K_{i}^{\ast}(t,s)(i=1,2)$ don't rely on
 $t$, we can easily infer that $D^{\alpha_{i}-1}T_{i}(u,v)(t)(i=1,2)$ is equiconvergent at
 $+\infty$. \par From the above three steps, Lemma 2.6 is satisfied. So the operator $T:P\rightarrow P$ is
relatively compact.
\par Finally we show that the operator $T:P\rightarrow P$ is continuous. Let $(u_{n},v_{n}),(u,v)\in P,$ such that $(u_{n},v_{n})\rightarrow
   (u,v)(n\rightarrow\infty).$ Then
   $||(u_{n},v_{n})||_{X\times Y}<+\infty,||(u,v)||_{X\times Y}<+\infty$. Similar to \eqref{y0} and \eqref{y00}, we
   have
\[\aligned
 ||T_{1}(u_{n},v_{n})||_{0}=&\sup_{t\in J}\Big|\int_{0}^{+\infty}\frac{K_{1}(t,s)}{1+t^{\alpha_{1}-1}}
 f_{1}(s,u_{n}(s),v_{n}(s),D^{\alpha_{1}-1}u_{n}(s),D^{\alpha_{2}-1}v_{n}(s))\mbox{d}s\Big|\\
 \leq&\frac{1}{\Gamma(\alpha_{1})-\Lambda_{1}}\Big[a_{10}^{\ast}+\sum_{k=1}^{4}a_{1k}^{\ast}||(u_{n},v_{n})||_{X\times
Y}^{\lambda_{1k}}\Big],
\endaligned \]
and
\[\aligned ||T_{1}(u_{n},v_{n})||_{1}=&\sup_{t\in
J}|\int_{0}^{+\infty}K^{\ast}_{1}(t,s)f_{1}(s,u_{n}(s),v_{n}(s),D^{\alpha_{1}-1}u_{n}(s),D^{\alpha_{2}-1}v_{n}(s))ds|\\
 \leq&\frac{\Gamma(\alpha_{1})}{\Gamma(\alpha_{1})-\Lambda_{1}}\Big[a_{10}^{\ast}+\sum_{k=1}^{4}a_{1k}^{\ast}||(u_{n},v_{n})||_{X\times
Y}^{\lambda_{1k}}\Big].
\endaligned \]
By continuity of function $f_{1}$ and the Lebesgue dominated
convergence theorem, we obtain
\[\aligned  &\lim_{n\rightarrow\infty}\int_{0}^{+\infty}\frac{K_{1}(t,s)}{1+t^{\alpha_{1}-1}}f_{1}(s,u_{n}(s),v_{n}(s),D^{\alpha_{1}-1}u_{n}(s),D^{\alpha_{2}-1}v_{n}(s))
\mbox{d}s\\=&\int_{0}^{+\infty}\frac{K_{1}(t,s)}{1+t^{\alpha_{1}-1}}f_{1}(s,u(s),v(s),D^{\alpha_{1}-1}u(s),D^{\alpha_{2}-1}v(s))\mbox{d}s,
\endaligned \]
and
\[\aligned&\lim_{n\rightarrow\infty}\int_{0}^{+\infty}K_{1}^{\ast}(t,s)f_{1}(s,u_{n}(s),v_{n}(s),D^{\alpha_{1}-1}u_{n}(s),D^{\alpha_{2}-1}v_{n}(s))\mbox{d}s
\\=&\int_{0}^{\infty}K_{1}^{\ast}(t,s)f_{1}(s,u(s),v(s),D^{\alpha_{1}-1}u(s),D^{\alpha_{2}-1}v(s))\mbox{d}s.
\endaligned \]
 Then
\[\aligned\|T_{1}(u_{n},v_{n})-T_{1}(u,v)\|_{0}
\leq&\sup_{t\in
J}\int_{0}^{+\infty}\frac{K_{1}(t,s)}{1+t^{\alpha-1}}\Big|f_{1}(s,u_{n}(s),v_{n}(s),D^{\alpha_{1}-1}u_{n}(s),D^{\alpha_{2}-1}v_{n}(s))
\\&-f_{1}(s,u(s),v(s),D^{\alpha_{1}-1}u(s),D^{\alpha_{2}-1}v(s))\Big|\mbox{d}s\rightarrow
0 ,\ n\rightarrow\infty ,\endaligned \] and
\[\aligned\|T_{1}(u_{n},v_{n})-T_{1}(u,v)\|_{1} \leq&\sup_{t\in J}\int_{0}^{+\infty}K^{\ast}_{1}(t,s)\Big|f_{1}(s,u_{n}(s),v_{n}(s),D^{\alpha_{1}-1}u_{n}(s),D^{\alpha_{2}-1}v_{n}(s))
\\&-f_{1}(s,u(s),v(s),D^{\alpha_{1}-1}u(s),D^{\alpha_{2}-1}v(s))\Big|\mbox{d}s\rightarrow
0 ,\ n\rightarrow\infty .\endaligned \] So, as $n\rightarrow\infty$,
\[\aligned
\|T_{1}(u_{n},v_{n})-T_{1}(u,v)\|_{X}=\max\{\|T_{1}(u_{n},v_{n})-T_{1}(u,v)\|_{0},
\|T_{1}(u_{n},v_{n})-T_{1}(u,v)\|_{1}\}\rightarrow 0.
\endaligned \]
 This means that the operator $T_{1}$ is continuous. At the same
 way, we can show than the operator $T_{2}$ is continuous. That is, the operator $T$ is
 continuous.
\par In view of the above all arguments, we deduce that the operator $T:P\rightarrow
P$ is completely continuous. Therefore proof is completed.

\section{Main results}

 For convenience, we set $$L_{i}=\displaystyle\frac{1}{\Gamma(\alpha_{i})-\Lambda_{i}}
,i=1,2,\ \
L=\displaystyle\max\{L_{1},L_{2},\Gamma(\alpha_{1})L_{1},\Gamma(\alpha_{2})L_{2}\}.$$

   Define a partial order over the product space:
 \begin{center}
   $ \begin{pmatrix}
                      u_{1}    \\
                     v_{1}    \\
   \end{pmatrix} $  $\geq$
  $ \begin{pmatrix}
                      u_{2}    \\
                     v_{2}    \\
   \end{pmatrix} $
   \end{center}
if $u_{1}(t)\geq u_{2}(t),v_{1}(t)\geq
v_{2}(t),D^{\alpha_{1}-1}u_{1}(t)\geq
D^{\alpha_{1}-1}u_{2}(t),D^{\alpha_{2}-1}v_{1}(t)\geq
D^{\alpha_{2}-1}v_{2}(t),t\in J.$\\
 {\bf Theorem 3.1.} \ Assume that
(H$_{1})$,(H$_{2})$ and (H$_{4})$ hold. There exists a positive
constant $R$ such that the system \eqref{fds} have two positive
solutions $(u^{\ast},v^{\ast})$ and $(w^{\ast},z^{\ast})$ satisfying
$0\leq \|(u^{\ast},v^{\ast})\|_{X\times Y}\leq R$ and $0\leq
\|(w^{\ast},z^{\ast})\|_{X\times Y}\leq R$. Moreover,
$\lim_{n\rightarrow\infty}(u_{n},v_{n})=(u^{\ast},v^{\ast})$ and
$\lim_{n\rightarrow\infty}(w_{n},z_{n})=(w^{\ast},z^{\ast})$,
 $(u_{n},v_{n})$ and $(w_{n},z_{n})$ can be given by the following
 monotone iterative schemes
\begin{equation}\label{t1}
 (u_{n},v_{n})=T(u_{n-1},v_{n-1})=
\left(\aligned
&{ T_{1}(u_{n-1},v_{n-1})(t)  } \\
&{  T_{2}(u_{n-1},v_{n-1})(t)}
\endaligned\right),n=1,2,\ldots,\
with \ (u_{0}(t),v_{0}(t))= \left(\aligned
&{  Rt^{\alpha_{1}} } \\
&{ R t^{\alpha_{2}}}
\endaligned\right)
 \end{equation}
 and
\begin{equation}\label{t2}
 (w_{n},z_{n})=T(w_{n-1},z_{n-1})=
\left(\aligned
&{ T_{1}(w_{n-1},z_{n-1})(t)  } \\
&{  T_{2}(w_{n-1},z_{n-1})(t)}
\endaligned\right),n=1,2,\ldots,\
with \ (w_{0}(t)z_{0}(t))= \left(\aligned
&{  0 } \\
&{0}
\endaligned\right).
 \end{equation}
 In addition
 $$\left(\aligned
&{w_{0}(t)} \\
&{ z_{0}(t)}
\endaligned\right)\leq
\left(\aligned
&{w_{1}(t)} \\
&{ z_{1}(t)}
\endaligned\right)\leq
\cdots\leq \left(\aligned
&{w_{n}(t)} \\
&{ z_{n}(t)}
\endaligned\right)
\leq\cdots\leq \left(\aligned
&{ w^{\ast}} \\
&{ z^{\ast}}
\endaligned\right)
\leq\cdots\leq \left(\aligned
&{  u^{\ast}} \\
&{ v^{\ast}}
\endaligned\right)
\leq\cdots\leq \left(\aligned
&{  u_{n}(t)} \\
&{  v_{n}(t)}
\endaligned\right)$$\begin{equation}\label{t3}
\leq\cdots\leq \left(\aligned
&{u_{2}(t)} \\
&{ v_{2}(t)}
\endaligned\right)
\leq \left(\aligned
&{u_{1}(t)} \\
&{ v_{1}(t)}
\endaligned\right)
 \end{equation}
and
 $$\left(\aligned
&{ D^{\alpha_{1}-1} w_{0}(t)} \\
&{ D^{\alpha_{2}-1} z_{0}(t)}
\endaligned\right)\leq
\left(\aligned
&{ D^{\alpha_{1}-1} w_{1}(t)} \\
&{ D^{\alpha_{2}-1} z_{1}(t)}
\endaligned\right)\leq\cdots\leq \left(\aligned
&{ D^{\alpha_{1}-1} w_{n}(t)} \\
&{ D^{\alpha_{2}-1} z_{n}(t)}
\endaligned\right)
\leq\cdots\leq \left(\aligned
&{ D^{\alpha_{1}-1} w^{\ast}} \\
&{ D^{\alpha_{2}-1}z^{\ast}}
\endaligned\right)
\leq\cdots\leq \left(\aligned
&{  D^{\alpha_{1}-1}u^{\ast}} \\
&{ D^{\alpha_{2}-1}v^{\ast}}
\endaligned\right)$$\begin{equation}\label{t4}
\leq\cdots\leq \left(\aligned
&{ D^{\alpha_{1}-1} u_{n}(t)} \\
&{ D^{\alpha_{2}-1} v_{n}(t)}
\endaligned\right)
\leq\cdots\leq \left(\aligned
&{ D^{\alpha_{1}-1} u_{2}(t)} \\
&{ D^{\alpha_{2}-1} v_{2}(t)}
\endaligned\right)
\leq \left(\aligned
&{ D^{\alpha_{1}-1} u_{1}(t)} \\
&{ D^{\alpha_{2}-1} v_{1}(t)}
\endaligned\right).
 \end{equation}
{\bf Proof.} \ First, Lemma 2.7 leads to the fact that $T(P) \subset
P$ for any $(u,v)\in P,t\in J$.
\par Next, for $0\leq \lambda_{1k},\lambda_{2k}<1(k=1,2,3,4)$, choose
$$ R\geq\max\Big\{5a_{10}^{\ast},5a_{20}^{\ast},(5La_{1k}^{\ast})^{1/(1-\lambda_{1k})},(5La_{2k}^{\ast})^{1/(1-\lambda_{2k})}\Big\},k=1,2,3,4,$$
and define $U_{R}=\{(u,v)\in P:||(u,v)||_{X\times Y}\leq R\}.$  For
any $(u,v)\in U_{R}$, similar to
 \eqref{y0} and \eqref{y00}, we obtain
$$
 ||T_{1}(u,v)||_{0}\leq
L_{1}\Big[a_{10}^{\ast}+\sum_{k=1}^{4}a_{1k}^{\ast}||(u,v)||_{X\times
Y}^{\lambda_{1k}}\Big]\leq
L\Big[a_{10}^{\ast}+\sum_{k=1}^{4}a_{1k}^{\ast}R^{\lambda_{1k}}\Big]\leq
R
 $$
and
$$||T_{1}(u,v)||_{1}\leq
L_{1}\Big[a_{10}^{\ast}+\sum_{k=1}^{4}a_{1k}^{\ast}||(u,v)||_{X\times
Y}^{\lambda_{1k}}\Big]\leq
L\Big[a_{10}^{\ast}+\sum_{k=1}^{4}a_{1k}^{\ast}R^{\lambda_{1k}}\Big]\leq
R.
 $$
 This implies that $||T_{1}(u,v)||_{X}\leq R$ for all $(u,v)\in
U_{R}$. In the same way, $||T_{2}(u,v)||_{Y}\leq R$.
 Consequently we
have
$$||T(u,v)||_{X\times Y}=\Big\{\|T_{1}(u,v)\|_{X},\|T_{2}(u,v)\|_{Y}\Big\}\leq R.$$
That is, $T(U_{R})\subset U_{R}$.

According to  \eqref{t1} and \eqref{t2}, it is obvious that
$(u_{0}(t),v_{0}(t)),(w_{0}(t),z_{0}(t))\in U_{R}$. By the complete
continuity of the operator $T$, we define the schemes
$(u_{n},v_{n})$ and $(w_{n},z_{n})$ by $(u_{n},v_{n})=
T(u_{n-1},v_{n-1}), (w_{n},z_{n})$= $T(w_{n-1},z_{n-1})$ for
$n=1,2,\ldots.$ Since $T(B)\subset B$, we can know that
$(u_{n},v_{n}),(w_{n},z_{n})\in T(B)$ for $n=1,2,\ldots.$ Hence we
need show that there exist $(u^{\ast},v^{\ast})$ and
$(w^{\ast},z^{\ast})$ satisfying
$\lim_{n\rightarrow\infty}(u_{n},v_{n})=(u^{\ast},v^{\ast})$ and
$\lim_{n\rightarrow\infty}(w_{n},z_{n})=(w^{\ast},z^{\ast})$, which
are two monotone schemes for positive solutions of the system
\eqref{fds}.

For $t\in J$, by Lemma 2.3 and \eqref{t1}, we know
\[\aligned
u_{1}(t)=T_{1}(u_{0},v_{0})(t)=&\displaystyle\int^{+\infty}_{0}K_{1}(t,s)f_{1}(s,u_{0}(s),v_{0}(s),
D^{\alpha_{1}-1}u_{0}(s),D^{\alpha_{2}-1}v_{0}(s))\mbox{d}s \\
\leq& \displaystyle t^{\alpha_{1}-1}
L_{1}\Big[a_{10}^{\ast}+\sum_{k=1}^{4}a_{1k}^{\ast}R^{\lambda_{1k}}\Big]\\
 \leq&Rt^{\alpha_{1}-1}=u_{0}(t)\endaligned \]
 and
\[\aligned
v_{1}(t)=T_{2}(u_{0},v_{0})(t)=&\displaystyle\int^{+\infty}_{0}K_{2}(t,s)f_{2}(s,u_{0}(s),v_{0}(s),
D^{\alpha_{1}-1}u_{0}(s),D^{\alpha_{2}-1}v_{0}(s))\mbox{d}s\\
\leq& \displaystyle t^{\alpha_{2}-1}
L_{2}\Big[a_{20}^{\ast}+\sum_{k=1}^{4}a_{2k}^{\ast}R^{\lambda_{2k}}\Big]\\
 \leq&Rt^{\alpha_{2}-1}=v_{0}(t),\endaligned \]
that is
\begin{equation}\label{pf1}
  T(u,v)(t)=
\left(\aligned
&{u_{1}(t)  } \\
&{  v_{1}(t)}
\endaligned\right)=
\left(\aligned
&{  T_{1}(u_{0},v_{0})(t) } \\
&{ T_{2}(u_{0},v_{0})(t) }
\endaligned\right)\leq
\left(\aligned
&{Rt^{\alpha_{1}-1}} \\
&{  Rt^{\alpha_{2}-1}}
\endaligned\right)=
\left(\aligned
&{ u_{0}(t)} \\
&{ v_{0}(t)}
\endaligned\right)
   \end{equation}

And then we study the monotonicity of the fractional derivative of
$(u,v)$. By \eqref{pf1} we know
\[\aligned
D^{\alpha_{1}-1}u_{1}(t)=&D^{\alpha_{1}-1}T_{1}(u_{0},v_{0})(t)
=\displaystyle\int^{+\infty}_{0}K^{\ast}_{1}(t,s)f_{1}(s,u_{0}(s),v_{0}(s)\
D^{\alpha_{1}-1}u_{0}(s),D^{\alpha_{2}-1}v_{0}(s))\mbox{d}s\\
\leq&
\Gamma(\alpha_{1})L_{1}\Big[a_{10}^{\ast}+\sum_{k=1}^{4}a_{1k}^{\ast}R^{\lambda_{1k}}\Big]
\leq\Gamma (\alpha_{1})R=D^{\alpha_{1}-1}u_{0}(t),\\
D^{\alpha_{2}-1}v_{1}(t)=&D^{\alpha_{2}-1}T_{2}(u_{0},v_{0})(t)
=\displaystyle\int^{+\infty}_{0}K^{\ast}_{2}(t,s)f_{2}(s,u_{0}(s),v_{0}(s)\
D^{\alpha_{1}-1}u_{0}(s),D^{\alpha_{2}-1}v_{0}(s))\mbox{d}s\\
\leq& \Gamma
(\alpha_{2})L_{2}\Big[a_{20}^{\ast}+\sum_{k=1}^{4}a_{2k}^{\ast}R^{\lambda_{1k}}\Big]
\leq\Gamma(\alpha_{2})R=D^{\alpha_{2}-1}v_{0}(t),\endaligned \]
 that
is
\begin{equation}\label{pf2}
  T(u,v)(t)=
\left(\aligned
&{D^{\alpha_{1}-1}u_{1}(t)  } \\
&{ D^{\alpha_{2}-1}v_{1}(t)}
\endaligned\right)=
\left(\aligned
&{  D^{\alpha_{1}-1}T_{1}(u_{0},v_{0})(t) } \\
&{ D^{\alpha_{2}-1}T_{2}(u_{0},v_{0})(t) }
\endaligned\right)\leq
\left(\aligned
&{ \Gamma(\alpha_{1})R } \\
&{ \Gamma(\alpha_{2})R}
\endaligned\right)=
\left(\aligned
&{D^{\alpha_{1}-1}u_{0}(t)} \\
&{  D^{\alpha_{2}-1}v_{0}(t)}
\endaligned\right)
   \end{equation}

Thus, from  \eqref{pf1} and  \eqref{pf2}, for $\forall t\in J$,  by
the monotonicity hypothesis (H$_{4})$ of the functions $f_{i}$, we
do the second iteration
\begin{center}
$\begin{pmatrix}
                      u_{2}(t)  \\
                     v_{2}(t) \\
   \end{pmatrix}  $=
$ \begin{pmatrix}
                      T_{1}(u_{1},v_{1})(t)   \\
                     T_{2}(u_{1},v_{1})(t)  \\
   \end{pmatrix}$
   $\leq$
$ \begin{pmatrix}
                      T_{1}(u_{0},v_{0})(t)   \\
                     T_{2}(u_{0},v_{0})(t)  \\
   \end{pmatrix}$
=
   $\begin{pmatrix}
                      u_{1}(t)  \\
                     v_{1}(t)  \\
   \end{pmatrix}  $,
   \end{center}

\begin{center}
$\begin{pmatrix}
                      D^{\alpha_{1}-1}u_{2}(t)  \\
                     D^{\alpha_{2}-1}v_{2}(t) \\
   \end{pmatrix}  $=
$ \begin{pmatrix}
                      D^{\alpha_{1}-1}T_{1}(u_{1},v_{1})(t)   \\
                     D^{\alpha_{2}-1}T_{2}(u_{1},v_{1})(t)  \\
   \end{pmatrix}$
   $\leq$
$ \begin{pmatrix}
                      D^{\alpha_{1}-1}T_{1}(u_{0},v_{0})(t)   \\
                     D^{\alpha_{2}-1}T_{2}(u_{0},v_{0})(t)  \\
   \end{pmatrix}$
=
   $\begin{pmatrix}
                      D^{\alpha_{1}-1}u_{1}(t)  \\
                     D^{\alpha_{2}-1}v_{1}(t)  \\
   \end{pmatrix}  $.
   \end{center}
By recursion, for $t\in J$, the scheme
$\{(u_{n},v_{n})\}_{n=0}^{\infty}$ satisfies
\begin{center}
$\begin{pmatrix}
                      u_{n+1}(t)  \\
                     v_{n+1}(t) \\
   \end{pmatrix}  $ $\leq$
$ \begin{pmatrix}
                      u_{n}(t)  \\
                     v_{n}(t) \\
   \end{pmatrix},$
   \ \
$ \begin{pmatrix}
                     D^{\alpha_{1}-1}u_{n+1}(t)  \\
                     D^{\alpha_{2}-1}v_{n+1}(t) \\
   \end{pmatrix}$
$\leq$
   $\begin{pmatrix}
                    D^{\alpha_{1}-1}u_{n}(t)  \\
                     D^{\alpha_{2}-1}v_{n}(t) \\
   \end{pmatrix}  $.
   \end{center}
By the aid of  the iterative scheme
$(u_{n+1},v_{n+1})=T(u_{n},v_{n})$ and the complete continuity of
the operator $T$, it is easy to infer that $(u_{n},v_{n})\rightarrow
(u^{\ast},v^{\ast})$ and $T(u^{\ast},v^{\ast})=(u^{\ast},v^{\ast})$
.

For the scheme $\{(w_{n},z_{n})\}_{n=0}^{\infty}$, we use a similar
discussion. For $t\in J$, we have
\begin{center}
$\begin{pmatrix}
                    w_{1}(t)  \\
                     z_{1}(t)  \\
   \end{pmatrix}  $=
$ \begin{pmatrix}
                    T_{1}(w_{0},z_{0})(t)  \\
                     T_{2}(w_{0},z_{0})(t) \\
   \end{pmatrix}$
= $ \begin{pmatrix}
                    \int^{+\infty}_{0}K_{1}(t,s)f_{1}(t,w_{0}(t),z_{0}(t),
D^{\alpha_{1}-1}w_{0}(t),D^{\alpha_{2}-1}z_{0}(t))ds   \\
                   \int^{+\infty}_{0}K_{2}(t,s)f_{2}(t,w_{0}(t),z_{0}(t),
D^{\alpha_{1}-1}w_{0}(t),D^{\alpha_{2}-1}z_{0}(t))ds \\
   \end{pmatrix}$\\
   $\geq$ $\begin{pmatrix}
                   0   \\
                    0  \\
   \end{pmatrix}  $=$\begin{pmatrix}
                      w_{0}(t)  \\
                     z_{0}(t)  \\
   \end{pmatrix}  $,~~~~~~~~~~~~~~~~~~
   \end{center}

\begin{center}
$\begin{pmatrix}
                    D^{\alpha_{1}-1}w_{1}(t)  \\
                     D^{\alpha_{2}-1}z_{1}(t)  \\
   \end{pmatrix}  $
   = $ \begin{pmatrix}
                    \int^{+\infty}_{0}K^{\ast}_{1}(t,s)f_{1}(t,w_{0}(t),z_{0}(t),
D^{\alpha_{1}-1}w_{0}(t),D^{\alpha_{2}-1}z_{0}(t))ds   \\
                   \int^{+\infty}_{0}K^{\ast}_{2}(t,s)f_{2}(t,w_{0}(t),z_{0}(t),
D^{\alpha_{1}-1}w_{0}(t),D^{\alpha_{2}-1}z_{0}(t))ds \\
   \end{pmatrix}$

   $\geq$ $\begin{pmatrix}
                   0   \\
                    0  \\
   \end{pmatrix}  $=$\begin{pmatrix}
                      D^{\alpha_{1}-1}w_{0}(t)  \\
                     D^{\alpha_{2}-1}z_{0}(t)  \\
   \end{pmatrix}  $.~~~~~~~~~~~~
   \end{center}
Using the the monotonicity hypothesis (H$_{4})$ of the functions
$f_{i}$, we have
\begin{center}
$\begin{pmatrix}
                    w_{2}(t)  \\
                     z_{2}(t)  \\
   \end{pmatrix}  $=
$ \begin{pmatrix}
                    T_{1}(w_{1},z_{1})(t)  \\
                     T_{2}(w_{1},z_{1})(t) \\
   \end{pmatrix}$
$\geq$ $ \begin{pmatrix}
                    T_{1}(w_{0},z_{0})(t)  \\
                     T_{2}(w_{0},z_{0})(t) \\
   \end{pmatrix}$
=$\begin{pmatrix}
                      w_{1}(t)  \\
                     z_{1}(t)  \\
   \end{pmatrix}  $,
   \end{center}

\begin{center}
$\begin{pmatrix}
                    D^{\alpha_{1}-1}w_{2}(t)  \\
                     D^{\alpha_{2}-1}z_{2}(t)  \\
   \end{pmatrix}  $=
$ \begin{pmatrix}
                    D^{\alpha_{1}-1}T_{1}(w_{1},z_{1})(t)  \\
                    D^{\alpha_{2}-1} T_{2}(w_{1},z_{1})(t) \\
   \end{pmatrix}$
$\geq$ $ \begin{pmatrix}
                    D^{\alpha_{1}-1}T_{1}(w_{0},z_{0})(t)  \\
                     D^{\alpha_{2}-1}T_{2}(w_{0},z_{0})(t) \\
   \end{pmatrix}$
=$\begin{pmatrix}
                      D^{\alpha_{1}-1}w_{1}(t)  \\
                     D^{\alpha_{2}-1}z_{1}(t)  \\
   \end{pmatrix}  $.
   \end{center}
Analogously, for $n=0,1,2,\ldots$ and $t\in J$, we have
\begin{center}
$\begin{pmatrix}
                      w_{n+1}(t)  \\
                     z_{n+1}(t) \\
   \end{pmatrix}  $ $\geq$
$ \begin{pmatrix}
                      w_{n}(t)  \\
                     z_{n}(t) \\
   \end{pmatrix},$
   \ \
$ \begin{pmatrix}
                     D^{\alpha_{1}-1}w_{n+1}(t)  \\
                     D^{\alpha_{2}-1}z_{n+1}(t) \\
   \end{pmatrix}$
$\geq$
   $\begin{pmatrix}
                    D^{\alpha_{1}-1}w_{n}(t)  \\
                     D^{\alpha_{2}-1}z_{n}(t) \\
   \end{pmatrix}  $.
   \end{center}
Combining the iterative scheme $(w_{n+1},z_{n+1})=T(w_{n},z_{n})$
and the complete continuity of the operator $T$, it is easy to infer
that $(w_{n},z_{n})\rightarrow (w^{\ast},z^{\ast})$ and
$T(w^{\ast},z^{\ast})=(w^{\ast},z^{\ast})$ . Finally we show that
$(u^{\ast},v^{\ast})$ and $(w^{\ast},z^{\ast})$ are the minimal and
maximal  positive solutions of the system  \eqref{fds}. Suppose that
$(\xi(t),\eta(t))$ is any positive solution of the system
\eqref{fds}, then $T(\xi(t),\eta(t))=(\xi(t),\eta(t))$ and
\begin{center}
$\begin{pmatrix}
                       w_{0}(t)  \\
                     z_{0}(t) \\
   \end{pmatrix}  $ =
   $\begin{pmatrix}
                      0  \\
                     0 \\
   \end{pmatrix}  $ $\leq$
$ \begin{pmatrix}
                   \xi(t)  \\
                   \eta(t) \\
   \end{pmatrix}$
$\leq$
   $\begin{pmatrix}
                   Rt^{\alpha_{1}-1}  \\
                    Rt^{\alpha_{2}-1}\\
   \end{pmatrix}  $
   =
   $\begin{pmatrix}
                  u_{0}(t)  \\
                    v_{0}(t)\\
   \end{pmatrix}  $,
   \end{center}

\begin{center}
$\begin{pmatrix}
                      D^{\alpha_{1}-1}w_{0}(t)  \\
                     D^{\alpha_{2}-1}z_{0}(t) \\
   \end{pmatrix}  $  $\leq$
$ \begin{pmatrix}
                   D^{\alpha_{1}-1}\xi(t) \\
                   D^{\alpha_{2}-1}\eta(t) \\
   \end{pmatrix}$
$\leq$
   $\begin{pmatrix}
                  D^{\alpha_{1}-1}u_{0}(t) \\
                    D^{\alpha_{2}-1}v_{0}(t)\\
   \end{pmatrix}  $.
   \end{center}
 Applying the monotone property of the operator $T$, we know that

\begin{center}
$ \begin{pmatrix}
                     w_{1}(t)  \\
                     z_{1}(t) \\
   \end{pmatrix}$ =$\begin{pmatrix}
                      T_{1}(w_{0},z_{0})(t)  \\
                    T_{2}(w_{0},z_{0})(t) \\
   \end{pmatrix}  $
$\leq$ $\begin{pmatrix}
                      \xi(t)  \\
                     \eta(t) \\
   \end{pmatrix}  $
 $\leq$ $\begin{pmatrix}
                      T_{1}(u_{0},v_{0})(t)  \\
                    T_{2}(u_{0},v_{0})(t) \\
   \end{pmatrix}  $
 =$\begin{pmatrix}
                      u_{1}(t)  \\
                    v_{1}(t) \\
   \end{pmatrix}  $,
   \end{center}
   \begin{center}
$ \begin{pmatrix}
                    D^{\alpha_{1}-1} w_{1}(t)  \\
                      D^{\alpha_{2}-1}z_{1}(t) \\
   \end{pmatrix}$
$\leq$ $\begin{pmatrix}
                     D^{\alpha_{1}-1} \xi(t)  \\
                      D^{\alpha_{2}-1}\eta(t) \\
   \end{pmatrix}  $
 $\leq$ $\begin{pmatrix}
                     D^{\alpha_{1}-1} u_{1}(t)  \\
                     D^{\alpha_{2}-1} v_{1}(t) \\
   \end{pmatrix}  $.
   \end{center}
Repeating the above steps, we have
\begin{center}
$ \begin{pmatrix}
                     w_{n}(t)  \\
                     z_{n}(t) \\
   \end{pmatrix}$
$\leq$ $\begin{pmatrix}
                      \xi(t)  \\
                     \eta(t) \\
   \end{pmatrix}  $
 $\leq$ $\begin{pmatrix}
                      u_{n}(t)  \\
                    v_{n}(t) \\
   \end{pmatrix}  $
   \end{center}
\begin{center}
$\begin{pmatrix}
                      D^{\alpha_{1}-1}w_{n}(t)  \\
                     D^{\alpha_{2}-1}z_{n}(t) \\
   \end{pmatrix}  $  $\leq$
$ \begin{pmatrix}
                   D^{\alpha_{1}-1}\xi(t) \\
                   D^{\alpha_{2}-1}\eta(t) \\
   \end{pmatrix}$
$\leq$
   $\begin{pmatrix}
                  D^{\alpha_{1}-1}u_{n}(t) \\
                    D^{\alpha_{2}-1}v_{n}(t)\\
   \end{pmatrix}  $,
   \end{center}
From the above results, combine $\lim\limits_{n\rightarrow\infty}
 (w_{n},z_{n})=(w^{\ast},z^{\ast})$ and
$\lim\limits_{n\rightarrow\infty}
(u_{n},u_{n})=(u^{\ast},v^{\ast})$, we get the results \eqref{t3}
and \eqref{t4}.

Again $f(t,0,0,0,0)\neq 0$ for all $t\in J$, we know that $(0,0)$
isn't a solution of the system (1). By \eqref{t3} and \eqref{t4}, it
is obvious
 that $(w^{\ast},z^{\ast})$ and $(u^{\ast},v^{\ast})$ are the
extreme positive solutions of system \eqref{fds}, which can be
constructed by means of two monotone iterative schemes in \eqref{t1}
and \eqref{t2}.
\par With regard to the difference scope of parameters $\lambda_{ik}(i=1,2,k=1,2,3,4)$, the method is similar, so we
omit the details, thus the proof is completed.

 {\bf Theorem 3.2.} \ Suppose the hypotheses (H$_{1})$ and (H$_{3})$ are
satisfied.
  If \begin{equation}\label{tt1} m=L\max\Big\{\sum_{k=1}^{4}b_{1k},\sum_{k=1}^{4}b_{2k}\Big\}<1, \end{equation} then
  the system \eqref{fds} has a unique positive solution
$(\overline{u}(t),\overline{v}(t))$ in $P$. Moreover, there is a
iterative scheme $(u_{n},v_{n})$, such that
$(u_{n},v_{n})\rightarrow (\overline{u},\overline{u})$ as
$n\rightarrow \infty$ uniformly on any finite interval of $J$, where

\begin{equation}\label{tt2}
(u_{n},v_{n})=T(u_{n-1},v_{n-1})=\left(\aligned
&{T_{1}(u_{n-1},v_{n-1})(t)} \\
&{T_{2}(u_{n-1},v_{n-1})(t)}
\endaligned\right),n=1,2,\ldots.
\end{equation}

In addition, there is an error estimate for the approximation
scheme.
\begin{equation}\label{tt3}||(u_{n},v_{n})-(\overline{u},\overline{v})||_{X\times Y}
=\frac{m^{n}}{1-m}||(u_{1},v_{1})-(u_{0},v_{0})||_{X\times Y},
n=1,2,\ldots.\end{equation}

{\bf Proof } \ Choose $$r\geq L\tau/(1-m),$$ where $m$ is defined by
\eqref{tt1} and $\tau=\{\tau_{1},\tau_{2}\},\tau_{i}$ is defined by
the hypothesis (H$_{3})$.

First we prove that $TU_{r}\subset U_{r}$, where $U_{r}=\{(u,v)\in
P,||(u,v)||_{X\times Y}\leq r\}$. For any $(u,v)\in U_{r}$, by Lemma
2.3, Remark 2.2 and Lemma 2.5, we have
\[\aligned
 ||T_{1}(u,v)||_{0}\leq&
 L\Big(\sum_{k=1}^{4}b_{1k}^{\ast}r+\tau_{1}\Big)
\endaligned \]
and
\[\aligned
||T_{1}(u,v)||_{1}
 \leq & L\Big(\sum_{k=1}^{4}b_{1k}^{\ast}r+\tau_{1}\Big),
\endaligned \]
which implies
$$||T_{1}(u,v)||_{X}\leq  L\Big(\sum_{k=1}^{4}b_{1k}^{\ast}r+\tau_{i}\Big)\leq mr+L\tau_{1}, \ \forall (u,v)\in U_{r}.$$
Similar
$$||T_{2}(u,v)||_{Y}\leq L\Big(\sum_{k=1}^{4}b_{2k}^{\ast}r+\tau_{2}\Big)\leq mr+L\tau_{2}, \ \forall (u,v)\in U_{r}.$$
So we have
$$||T(u,v)||_{X\times Y}\leq mr+L\tau\leq r. \ \forall (u,v)\in U_{r}.$$

Now we show that T is a contraction. For any
$(u_{1},v_{1}),(u_{2},v_{2})\in U_{r}$, by hypothesis (H$_{3})$, we
obtain
\[\aligned
& ||T_{1}(u_{1},v_{1})-T_{1}(u_{2},v_{2})||_{0}\\ \leq &\sup_{t\in
 J}\int_{0}^{+\infty}\frac{K_{1}(t,s)}{1+t^{\alpha_{1}-1}}\Big|f_{1}(s,u_{1}(s),v_{1}(s),D^{\alpha_{1}-1}u_{1}(s),D^{\alpha_{2}-1}v_{1}(s))\\
 &-f_{1}(s,u_{2}(s),v_{2}(s),D^{\alpha_{1}-1}u_{2}(s),D^{\alpha_{2}-1}v_{2}(s))\Big|ds\\
\leq&
L\int_{0}^{+\infty}\Big[b_{11}(s)(1+s^{\alpha_{1}-1})\frac{|u_{1}(s)-u_{2}(s)|}{1+s^{\alpha_{1}-1}}
+b_{12}(s)(1+s^{\alpha_{2}-1})\frac{|v_{1}(s)-v_{2}(s)|}{1+s^{\alpha_{2}-1}}\\
&+b_{13}(s)|D^{\alpha_{1}-1}u_{1}(s)-D^{\alpha_{1}-1}u_{2}(s)|\Big]ds
+b_{14}(s)|D^{\alpha_{2}-1}v_{1}(s)-D^{\alpha_{2}-1}v_{2}(s)|\Big]ds\\
\leq
&L\sum_{k=1}^{4}b^{\ast}_{1k}||(u_{1},v_{1})-(u_{2},v_{2})||_{X\times
Y}
\endaligned \]
and
\[\aligned
||T_{1}(u_{1},v_{1})-T_{1}(u_{2},v_{2})||_{1}\leq & \sup_{t\in
J}\int_{0}^{+\infty}K_{1}^{\ast}(t,s)\Big|f_{1}(s,u_{1}(s),v_{1}(s),D^{\alpha_{1}-1}u_{1}(s),D^{\alpha_{2}-1}v_{1}(s))\\
 &-f_{1}(s,u_{2}(s),v_{2}(s),D^{\alpha_{1}-1}u_{2}(s),D^{\alpha_{2}-1}v_{2}(s))\Big|ds\\
 \leq
&L\sum_{k=1}^{4}b^{\ast}_{1k}||(u_{1},v_{1})-(u_{2},v_{2})||_{X\times
Y},
\endaligned \]
which implies
\begin{equation}\label{pp1}||T_{1}(u_{1},v_{1})-T_{1}(u_{2},v_{2})||_{X}\leq L\sum_{k=1}^{4}b^{\ast}_{1k}||(u_{1},v_{1})-(u_{2},v_{2})||_{X\times
Y}.\end{equation} In the same way, we have
\begin{equation}\label{pp2}||T_{2}(u_{1},v_{1})-T_{2}(u_{2},v_{2})||_{Y}\leq L\sum_{k=2}^{4}b^{\ast}_{2k}||(u_{1},v_{1})-(u_{2},v_{2})||_{X\times
Y}.\end{equation} From \eqref{pp1} and \eqref{pp2}, we have
\begin{equation}\label{pp3}||T(u_{1},v_{1})-T(u_{2},v_{2})||_{X\times Y}\leq m||(u_{1},v_{1})-(u_{2},v_{2})||_{X\times
Y},,\forall (u_{1},v_{1}),(u_{2},v_{2})\in U_{r}.\end{equation}
Since $m<1$, then T is a contraction.  Hence the Banach fixed-point
theorem ensures that $T$ has a unique fixed point
$(\overline{u},\overline{v})$ in $P$. That is, the system (1) has a
unique positive solution $(\overline{u},\overline{v})$.
\par Furthermore, for any $(u_{0},v_{0})\in P,\|(u_{n},v_{n})
-(\overline{u},\overline{v})\|_{X\times Y}\rightarrow 0$ as
$n\rightarrow\infty$, where $u_{ n}=T_{1}(u_{n-1},v_{n-1}),v_{
n}=T_{2}(u_{n-1},V_{n-1}),n=1,2,\ldots.$ By \eqref{pp3}, we obtain
$$||(u_{n},v_{n})-(u_{n-1},v_{n-1})||_{X\times Y}\leq m^{n-1}||(u_{1},v_{1})-(u_{0},v_{0})||_{X\times Y},$$
and
\begin{equation}\label{pp4}
\begin{array}{rcl}
 ||(u_{n},v_{n})-(u_{j},v_{j})||_{X\times Y}&\leq & ||(u_{n},v_{n})-(u_{n-1},v_{n-1})||_{X\times Y}+
 ||(u_{n-1},v_{n-1})-(u_{n-2},v_{n-2})||_{X\times Y}\\&&+\cdots + ||(u_{j+1},v_{j+1})-(u_{j},v_{j})||_{X\times Y}\\
&\leq
&\displaystyle\frac{m^{n}(1-m^{j-n})}{1-m}||(u_{1},v_{1})-(u_{0},v_{0})||_{X\times
Y}.
\end{array}\end{equation}
Letting $j\rightarrow +\infty$ on both sides of \eqref{pp4}, we have
$$||(u_{n},v_{n})-(\overline{u},\overline{v})||_{X\times Y}\leq \frac{m^{n}}{1-m}||u_{1}-u_{0}||_{X\times Y}.
$$
Hence the proof of theorem 3.2 is completed.

Now we give two examples to illustrate the application of the main
results.

 {\bf Example 3.1.} \ Consider the following fractional
differential system on an infinite interval
 \begin{equation}\label{ex1} \left\{
      \begin{array}{ll}
       \displaystyle
       -D^{2.5}u(t)=\frac{2}{(10+t)^{2}}+\frac{e^{-t}|u(t)|^{0.1}}{(1+\sqrt[]{t^{3}})^{0.1}}
       +\frac{e^{-2t}|v(t)|^{0.3}}{(1+\sqrt[]{t})^{0.3}}
       +\frac{2t| D^{1.5}u(t)|^{0.2}}{(3+t^{2})^{2}}
        +\frac{| D^{0.5}v(t)|^{0.4}}{1+t^{2}},\\
        -\displaystyle D^{1.5}v(t)=\frac{1}{(20+t)^{3}}+\frac{e^{-3t}|u(t)|^{0.2}}{(1+\sqrt[]{t^{3}})^{0.2}}
       +\frac{e^{-4t}|v(t)|^{0.4}}{(1+\sqrt[]{t})^{0.4}}
       +\frac{3t^{2}| D^{1.5}u(t)|^{0.2}}{(3+t^{3})^{2}}
        +\frac{2| D^{0.5}v(t)|^{0.6}}{1+t^{2}},\\
       u(0)=u'(0)=0,\ D^{1.5}u(+\infty)=
       \displaystyle\int_{0}^{+\infty}t^{-1.5}e^{-t}u(t)dt,\\
       u(0)=0,\ D^{0.5}v(+\infty)=
       \displaystyle\int_{0}^{+\infty}t^{-0.5}e^{-2t}v(t)dt,
     \end{array}
   \right.\end{equation}
   where $\alpha_{1}=2.5,\alpha_{1}=1.5,h_{1}(t)=t^{-1.5}e^{-t},h_{2}(t)=t^{-0.5}e^{-2t},\lambda_{11}=0.1,\lambda_{12}=0.3,
   \lambda_{13}=0.2,\lambda_{14}=0.4,\lambda_{21}=0.2,\lambda_{22}=0.4,\lambda_{23}=0.2,\lambda_{24}=0.6$ and
$$
f_{1}(t,u_{1},u_{2},u_{3},u_{4})=\frac{2}{(10+t)^{2}}+\frac{e^{-t}|u_{1}|^{0.1}}{(1+\sqrt[]{t^{3}})^{0.1}}
       +\frac{e^{-2t}|u_{2}|^{0.3}}{(1+\sqrt[]{t})^{0.3}}
       +\frac{2t|u_{3}|^{0.2}}{(3+t^{2})^{2}}
        +\frac{| u_{4}|^{0.4}}{1+t^{2}},$$
     $$ f_{2}(t,u_{1},u_{2},u_{3},u_{4})=\frac{1}{(20+t)^{3}}+\frac{e^{-3t}|u_{1}|^{0.2}}{(1+\sqrt[]{t^{3}})^{0.2}}
       +\frac{e^{-4t}|u_{2}|^{0.4}}{(1+\sqrt[]{t})^{0.4}}
       +\frac{3t^{2}| u_{3}|^{0.2}}{(3+t^{3})^{2}}
        +\frac{2| u_{4})|^{0.6}}{1+t^{2}},$$
   \par It is easy to know that $\Gamma
   (2.5)=1.32934>\Lambda_{1}=\int_{0}^{+\infty}h_{1}(t)t^{1.5}dt=1,\  \Gamma
   (1.5)=0.88623>\Lambda_{2}=\int_{0}^{+\infty}h_{2}(t)t^{0.5}dt=0.5,f_{i}(t,0,0,0,0)\not\equiv 0,i=1,2$.
   So the hypothesis
   (H$_{1})$ is satisfied.
   \par Noting that
$$\begin{array}{rcl}
|f_{1}(t,u_{1},u_{2},u_{3},u_{4})|&\leq&\displaystyle\frac{2}{(10+t)^{2}}+\frac{e^{-t}|u_{1}|^{0.1}}{(1+\sqrt[]{t^{3}})^{0.1}}
       +\frac{e^{-2t}|u_{2}|^{0.3}}{(1+\sqrt[]{t})^{0.3}}
       +\frac{2t|u_{3}|^{0.2}}{(3+t^{2})^{2}}
        +\frac{| u_{4}|^{0.4}}{1+t^{2}}\\
&=&a_{10}(t)+a_{11}(t)|u_{1}|^{0.1}+a_{12}(t)|u_{2}|^{0.3}+a_{13}(t)|u_{3}|^{0.2}+a_{14}(t)|u_{4}|^{0.4},
\end{array}$$
    $$\begin{array}{rcl}|f_{2}(t,u_{1},u_{2},u_{3},u_{4})|&\leq&\displaystyle\frac{1}{(20+t)^{3}}
    +\frac{e^{-3t}|u_{1}|^{0.2}}{(1+\sqrt[]{t^{3}})^{0.2}}
       +\frac{e^{-4t}|u_{2}|^{0.4}}{(1+\sqrt[]{t})^{0.4}}
       +\frac{3t^{2}| u_{3}|^{0.2}}{(3+t^{3})^{2}}
        +\frac{2| u_{4})|^{0.6}}{1+t^{2}}\\
&{\color{red}=}&a_{20}(t)+a_{21}(t)|u_{1}|^{0.2}+a_{22}(t)|u_{2}|^{0.2}+a_{23}(t)|u_{3}|^{0.2}+a_{24}(t)|u_{4}|^{0.6}
\end{array}$$
and
$$a_{10}^{\ast}=\displaystyle\int_{0}^{+\infty}a_{10}(t)dt=\frac{1}{5},\
 a_{11}^{\ast}=\displaystyle\int_{0}^{+\infty}a_{11}(t)(1+t^{1.5})^{0.1}dt=1,\
  a_{12}^{\ast}=\displaystyle\int_{0}^{+\infty}a_{12}(t)(1+t^{0.5})^{0.3}dt=\frac{1}{2},$$
$$a_{13}^{\ast}=\displaystyle\int_{0}^{+\infty}a_{13}(t)dt=\frac{1}{3},\
a_{14}^{\ast}=\displaystyle\int_{0}^{+\infty}a_{14}(t)dt=\frac{\pi}{2},$$
$$a_{20}^{\ast}=\displaystyle\int_{0}^{+\infty}a_{10}(t)dt=\frac{1}{800},\
 a_{21}^{\ast}=\displaystyle\int_{0}^{+\infty}a_{21}(t)(1+t^{1.5})^{0.2}dt=\frac{1}{3},\
  a_{22}^{\ast}=\displaystyle\int_{0}^{+\infty}a_{22}(t)(1+t^{0.5})^{0.4}dt=\frac{1}{4},$$
$$a_{23}^{\ast}=\displaystyle\int_{0}^{+\infty}a_{23}(t)dt=\frac{1}{3},\
a_{24}^{\ast}=\displaystyle\int_{0}^{+\infty}a_{24}(t)dt=\pi.$$
which means that the hypothesis (H$_{2})$ is satisfied.
\par From the expression of the function $f_{i}$, it is obvious that $f_{i}$ is increasing respect to the variables
$u_{1},u_{2},u_{3},u_{4},\forall t\in J, i=1,2$. Thus the hypothesis
(H$_{4})$ is satisfied. By Theorem 3.1, it follows that the system
\eqref{ex1} have two positive solution, which can be given by the
limits means of two explicit monotone iterative scheme in \eqref{t1}
and \eqref{t2}.

{\bf Example 3.2.} \ Consider the following fractional differential
system an infinite interval
\begin{equation}\label{ex2} \left\{
      \begin{array}{ll}
       \displaystyle
       -D^{2.5}u(t)=\frac{2}{(10+t)^{2}}+\frac{e^{-20t}|u(t)|}{1+\sqrt[]{t^{3}}}
       +\frac{e^{-15t}|v(t)|}{1+\sqrt[]{t}}
       +\frac{t| D^{1.5}u(t)|}{5(3+t^{2})^{2}}
        +\frac{t| D^{0.5}v(t)|}{10(1+t^{2})^{2}},\\
        -\displaystyle D^{1.5}v(t)=\frac{1}{(20+t)^{3}}+\frac{e^{-18t}|u(t)|}{1+\sqrt[]{t^{3}}}
       +\frac{e^{-16t}|v(t)|}{1+\sqrt[]{t}}
       +\frac{3t^{2}| D^{1.5}u(t)|}{7(3+t^{3})^{2}}
        +\frac{| D^{0.5}v(t)|}{20(1+t^{2})^{2}},\\
       u(0)=u'(0)=0,\ D^{1.5}u(+\infty)=
       \displaystyle\int_{0}^{+\infty}t^{-1.5}e^{-t}u(t)dt,\\
       u(0)=0,\ D^{0.5}v(+\infty)=
       \displaystyle\int_{0}^{+\infty}t^{-0.5}e^{-2t}v(t)dt,
     \end{array}
   \right.\end{equation}
    where $\alpha_{1}=2.5,\alpha_{1}=1.5,h_{1}(t)=t^{-1.5}e^{-t},h_{2}(t)=t^{-0.5}e^{-2t}$ and
  $$
f_{1}(t,u_{1},u_{2},u_{3},u_{4})=\frac{2}{(10+t)^{2}}+\frac{e^{-20t}|u_{1}|}{1+\sqrt[]{t^{3}}}
       +\frac{e^{-15t}|u_{2}|}{1+\sqrt[]{t}}
       +\frac{t|u_{3}|}{5(3+t^{2})^{2}}
        +\frac{t| u_{4}|}{10(1+t^{2})^{2}},$$
     $$ f_{2}(t,u_{1},u_{2},u_{3},u_{4})=\frac{1}{(20+t)^{3}}+\frac{e^{-18t}|u_{1}|}{1+\sqrt[]{t^{3}}}
       +\frac{e^{-16t}|u_{2}|}{1+\sqrt[]{t}}
       +\frac{3t^{2}| u_{3}|}{7(3+t^{3})^{2}}
        +\frac{| u_{4})|}{20(1+t^{2})},$$
   \par Similar to the example 3.1, it is  easy to verify that the hypothesis
   (H$_{1})$ is satisfied.
   \par Observing that
$$\begin{array}{rcl}
&&|f_{1}(t,u_{1},u_{2},u_{3},u_{4})-f_{1}(t,\overline{u}_{1},\overline{u}_{2},\overline{u}_{3},\overline{u}_{4})|\\
&\leq&\displaystyle\frac{e^{-20t}}{1+\sqrt[]{t^{3}}}|u_{1}-\overline{u}_{1}|
       +\frac{e^{-15t}|}{1+\sqrt[]{t}}|u_{2}-\overline{u}_{2}|+\frac{t}{5(3+t^{2})^{2}}|u_{3}-\overline{u}_{3}|
       +\frac{t}{10(1+t^{2})^{2}}|u_{4}-\overline{u}_{4}|\\
&=
&b_{11}(t)|u_{1}-\overline{u}_{1}|+b_{12}(t)|u_{2}-\overline{u}_{2}|+b_{13}(t)|u_{3}-\overline{u}_{3}|
+b_{14}(t)|u_{4}-\overline{u}_{4}|,
\end{array}$$
$$\begin{array}{rcl}
&&|f_{2}(t,u_{1},u_{2},u_{3},u_{4})-f_{2}(t,\overline{u}_{1},\overline{u}_{2},\overline{u}_{3},\overline{u}_{4})|\\
&\leq&\displaystyle\frac{e^{-18t}}{1+\sqrt[]{t^{3}}}|u_{1}-\overline{u}_{1}|
       +\frac{e^{-16t}|}{1+\sqrt[]{t}}|u_{2}-\overline{u}_{2}|+\frac{3t^{2}}{7(3+t^{3})^{2}}|u_{3}-\overline{u}_{3}|
       +\frac{1}{20(1+t^{2})}|u_{4}-\overline{u}_{4}|\\
&=&b_{21}(t)|u_{1}-\overline{u}_{1}|+b_{22}(t)|u_{2}-\overline{u}_{2}|+b_{23}(t)|u_{3}-\overline{u}_{3}|
+b_{24}(t)|u_{4}-\overline{u}_{4}|,
\end{array}$$
and
$$b_{11}^{\ast}=\displaystyle\int_{0}^{+\infty}b_{11}(t)(1+t^{1.5})dt=\frac{1}{20},\
  b_{12}^{\ast}=\displaystyle\int_{0}^{+\infty}b_{12}(t)(1+t^{0.5})dt=\frac{1}{15},$$
$$b_{13}^{\ast}=\displaystyle\int_{0}^{+\infty}b_{13}(t)dt=\frac{1}{30},\
b_{14}^{\ast}=\displaystyle\int_{0}^{+\infty}b_{14}(t)dt=\frac{1}{20},$$
$$b_{21}^{\ast}=\displaystyle\int_{0}^{+\infty}b_{21}(t)(1+t^{1.5})dt=\frac{1}{18},\
  b_{22}^{\ast}=\displaystyle\int_{0}^{+\infty}b_{22}(t)(1+t^{0.5})dt=\frac{1}{16},$$
$$b_{23}^{\ast}=\displaystyle\int_{0}^{+\infty}b_{23}(t)dt=\frac{1}{21},\
b_{24}^{\ast}=\displaystyle\int_{0}^{+\infty}a_{24}(t)dt=\frac{\pi}{40}.$$
$$\lambda_{1}=\displaystyle\int_{0}^{+\infty}f_{1}(t,0,0,0,0)dt=\displaystyle\int_{0}^{+\infty}\frac{2}{(10+t)^{2}}dt=\frac{1}{5},$$
$$\lambda_{2}=\displaystyle\int_{0}^{+\infty}f_{2}(t,0,0,0,0)dt=\displaystyle\int_{0}^{+\infty}\frac{1}{(20+t)^{3}}dt=\frac{\pi}{8000},$$
which means that the hypothesis (H$_{3})$ is satisfied. By direct
computation, we have
$$ m=L\max\Big\{\sum_{k=1}^{4}b_{1k},\sum_{k=1}^{4}b_{2k}\Big\}=4.03638\times
\max\Big\{0.2,0.24422\Big\}=0.98576<1.$$ So all conditions of
Theorem 3.2 are satisfied. Then the system \eqref{ex2} has a unique
positive solution, which can be obtained by the limits from the
iterative sequences in \eqref{tt2}.

\section{Conclusions}

In this paper, we apply the monotone iterative technique and the
Banach contraction mapping principle to study a class of fractional
differential system with integral boundary in an infinite interval.
We first transform the system \eqref{fds} into an equivalent
operator equation \eqref{op}, and then construct some norm
inequalities related to nonlinear terms $f_{i}(i=1,2)$ by means of
hypothesis conditions. Finally some explicit monotone iterative
schemes for approximating the extreme positive solutions and the
unique positive solution are established.

\vskip.20in

\noindent{\bf Acknowledgements}

 This work is supported by University
Natural Science Foundation of Anhui Provincial Education Department
(Grant No. KJ2018A0452), the Foundation of Suzhou University (Grant
No. 2019XJZY02, 2019XJSN03), Technology Research Foundation of
Chongqing Educational Committee(Grant No. KJQN201800533).

\noindent{\bf Conflicts of Interest}

The authors declare no conflict of
interest.\\

\end{document}